\documentclass[10pt,notitlepage]{article}

\usepackage{amssymb}
\usepackage[intlimits]{amsmath}
\usepackage{amsfonts}
\usepackage{amsthm} 

\usepackage{mathptmx}

\usepackage{hyperref}

\usepackage{cite}

\usepackage[a4paper]{geometry}

\let\oldsqrt\sqrt
\def\sqrt{\mathpalette\DHLhksqrt}
\def\DHLhksqrt#1#2{%
\setbox0=\hbox{$#1\oldsqrt{#2\,}$}\dimen0=\ht0
\advance\dimen0-0.2\ht0
\setbox2=\hbox{\vrule height\ht0 depth -\dimen0}%
{\box0\lower0.4pt\box2}}

\allowdisplaybreaks

\newcommand{\R}{\mathbb{R}} % reelle Zahlen
\newcommand{\N}{\mathbb{N}} % natuerliche Zahlen
\newcommand{\dist}{\textnormal{dist}} % dist ...
\newcommand{\diam}{\textnormal{diam}} % diam ...
\newcommand{\supp}{\textnormal{supp}} % supp ...
\newcommand{\tr}{\textnormal{tr}} % Trace ...
\newcommand{\inrad}{\textnormal{inrad}} % inrad ...
\newcommand{\hl}{(-\Delta)^{s}} % alpha/2 Laplace
\newcommand{\dom}{\textnormal{dom}} % Domain
\newcommand{\osc}{\textnormal{osc}}

\renewcommand{\H}{{\mathcal H}}

\renewcommand{\phi}{\varphi}

\renewcommand{\div}{\textnormal{div}}

\newcommand{\cB}{{\mathcal B}}

\newcommand{\cE}{{\mathcal E}}

\newcommand{\cH}{{\mathcal H}}

\newcommand{\cK}{{\mathcal K}}

\newcommand{\cV}{{\mathcal V}}

\newcommand{\eps}{\varepsilon}

\newtheorem{theorem}{Theorem}[section]
\newtheorem{corollary}[theorem]{Corollary}

\newtheorem{lemma}[theorem]{Lemma}
\newtheorem{proposition}[theorem]{Proposition}

\theoremstyle{definition}
\newtheorem{definition}[theorem]{Definition}
\newtheorem{remark}[theorem]{Remark}

\title{Asymptotic symmetry for a class of nonlinear fractional
  reaction-diffusion equations}
\author{Sven Jarohs\footnote{Institut f\"ur Mathematik, Goethe-Universit\"at, Frankfurt, Robert-Mayer-Stra\ss e 10, D-60054 Frankfurt, jarohs@math.uni-frankfurt.de.}
\hspace{1ex}and Tobias Weth\footnote{Institut f\"ur Mathematik, Goethe-Universit\"at, Frankfurt, Robert-Mayer-Stra\ss e 10, D-60054 Frankfurt, weth@math.uni-frankfurt.de.}
}

\begin{document}
\maketitle

\begin{abstract}
We study the
nonlinear fractional reaction-diffusion equation $\partial_{t}u + \hl u=
f(t,x,u)$, $s\in(0,1)$ in a bounded domain $\Omega$ together
with Dirichlet boundary conditions on \mbox{$\R^N \setminus \Omega$.} We
prove asymptotic symmetry of nonnegative globally bounded solutions
in the case where the underlying data obeys some symmetry and
monotonicity assumptions. More precisely, we assume that $\Omega$ is symmetric with respect to reflection at a hyperplane, say $\{x_1=0\}$, 
and convex in the $x_1$-direction, and that the
nonlinearity $f$ is even in $x_1$ and nonincreasing in $|x_1|$. Under
rather weak additional technical assumptions, we then show that any nonzero 
element in the $\omega$-limit set of nonnegative globally bounded
solution is even in $x_1$ and strictly
decreasing in $|x_1|$. This result, which is obtained via a series of
new estimates for antisymmetric supersolutions of a corresponding
family of linear equations, implies a strong maximum type
principle which is not available in the non-fractional case $s=1$. 
\end{abstract}
{\footnotesize
\begin{center}
\end{center}
\begin{center}
\end{center}
}

\section{Introduction}
\label{sec:introduction}

We consider the nonlinear fractional diffusion problem
\begin{equation*}
 (P)\qquad
\left\{
 \begin{aligned}
 \partial_{t}u + \hl u&= f(t,x,u) &&\qquad \text{in $(0,\infty)
   \times \Omega $,} \\
             u&=0                  &&\qquad \text{ on  $(0,\infty)
   \times (\R^{N}\setminus \Omega),$}
            \  
 \end{aligned}
\right.
\end{equation*}
where $\Omega$ is a bounded Lipschitz domain in $\R^{N}$,
 $s\in(0,1)$ and $f$ is a nonlinearity defined on
$(0,\infty) \times \Omega\times \cB$. Here and in the following, $\cB \subset
\R$ is an open interval (further assumptions on $f$ are to be specified later). Moreover, $\hl$ denotes the fractional
Laplacian, which for functions $u \in H^{2s}(\R^N)$ is defined via
Fourier transform:
\begin{equation}
  \label{eq:1}
 \widehat {\hl u}(\xi)=|\xi|^{2s} \widehat u (\xi) \qquad \text{for a.e.
  $\xi \in \R^N$.} 
\end{equation}
Fueled by various
applications in physics, biology or finance, linear and nonlinear
equations of the form given in $(P)$ or of similar type have
received immensely growing attention recently. In particular, evolution equations involving the fractional Laplacian appear in
the quasi-geostrophic equations (see e.g. \cite{CV10,SS03}) and in the fractional porous medium equation (see
\cite{PQRV11}), while further applications in the context of stable
processes are considered e.g. in \cite{A09,J05}. Very recently, the fractional Laplacian has been
studied in conformal geometry, see \cite{chang-gonzales,JX11}. In order to incorporate the Dirichlet boundary
condition on $\R^N \setminus \Omega$ in $(P)$, the operator $\hl$
has to be replaced by the Friedrichs extension of the restriction of
$\hl$ to the space $C_c^\infty(\Omega) \subset L^2(\Omega)$ of test
functions. Here and in the following, we identify $L^2(\Omega)$ with the space $\{u \in L^2(\R^N)\::\: \text{$u \equiv 0$ on $\R^N \setminus \Omega$}\}$.
This new operator, which we will also denote by $\hl$ in the
following, 
has the form domain $\mathcal{H}^{s}_{0}(\Omega)= \{u \in
H^s(\R^N)\::\: \text{$u \equiv 0$ on $\R^N \setminus \Omega$}\}$, and
it is widely used in analysis and probability
theory. In particular, it has recently been considered in the context
of semilinear problems, see e.g. \cite{BLW05,SV,BCD} and the
references therein. In probabilistic terms, the operator coincides with the generator of the $2s$-stable process in
$\Omega$ killed upon leaving $\Omega$. We note that for $u\in
C_c^\infty(\Omega)$ we have the
representation 
\begin{equation}
  \label{eq:59}
 \hl u(x)= c_{N,s}\, P.V.\!\int_{\R^N}\frac{u(x)-u(y)}{|x-y|^{N+2s}}dy
=c_{N,s}\, P.V.\!\int_{\Omega}\frac{u(x)-u(y)}{|x-y|^{N+2s}}dy +
u(x)\kappa_{\Omega}(x)  
\end{equation}
for $x \in \Omega$, where $P.V.$ stands for the principal value
integral and
\begin{equation}
  \label{eq:45}
c_{N,s}=s(1-s)\pi^{-N/2}4^s\frac{\Gamma(\frac{N}{2}+s)}{\Gamma(2-s)},\qquad 
\kappa_{\Omega}(x):=c_{N,s}\int\limits_{\R^{N}\setminus\Omega}|x-y|^{-N-2s}dy
\; \text{for $x \in \Omega$,}
\end{equation}
see e.g. \cite[Remark 3.11]{CS}. 

The focus of the present paper is the asymptotic shape of
global bounded solutions of $(P)$, i.e., the symmetry (and monotonicity)
of elements in the corresponding $\omega$-limit sets. For this we will use a weak
formulation for solutions of $(P)$. 
The quadratic form corresponding to $(-\Delta)^s$ on
$\mathcal{H}^{s}_{0}(\Omega)$ is given by 
\begin{align}
\mathcal{E}(u,v)&=\langle (-\Delta)^{s/2}u,(-\Delta)^{s/2}v\rangle_{L^2(\R^N)}\notag\\
&=  \frac{c_{N,s}}{2}\int_{\R^{N}}\int_{\R^{N}}\frac{(u(x)-u(y))(v(x)-v(y))}{|x-y|^{N+2s}}\ dx\ dy.\label{eq:46}
\end{align}
Since $\Omega$ is bounded, $\mathcal{E}$ defines a scalar product on
$\mathcal{H}^{s}_{0}(\Omega)$ which is equivalent to the standard scalar product
induced from the embedding $\mathcal{H}^{s}_{0}(\Omega)
\hookrightarrow H^s(\R^N)$. Consider the space $C_0(\Omega):= \{u \in
C(\R^N)\::\: \text{$u \equiv 0$ on $\R^N \setminus \Omega$}\}$
endowed with the usual $L^\infty$-norm. We
say that a function $u: (0,\infty) \times \R^N \to \R$ is a solution of $(P)$ if $u \in
C((0,\infty),\mathcal{H}^{s}_0(\Omega)\cap C_{0}(\Omega))\cap
C^{1}((0,\infty),L^{2}(\Omega))$, $u(t,x) \in \cB$ for every $(t,x) \in (0,\infty) \times \Omega$
and  
\begin{equation}\label{ungleichung}
 \mathcal{E}(u(t),\varphi) = \int_{\Omega} (f(t,x,u)-\partial_{t}
u)\varphi\ dx \qquad\text{ for all } \phi \in
\mathcal{H}^{s}_{0}(\Omega),\: t\in (0,\infty).
\end{equation}
For a solution $u$ of
$(P)$, we define the $\omega$-limit set (with respect to the norm
$\|\cdot\|_{L^\infty}$) as
\[
 \omega(u):=\left\{ z \in C_0(\Omega):\|u(t_{k})-z\|_{L^\infty}\to0\text{ for some }t_{k}\to\infty\right\}
\]
To state our main result, we introduce the following assumptions.
\begin{enumerate}
\item[(D1)] $\Omega$ is bounded with a Lipschitz boundary. Moreover,
  $\Omega$ is convex and symmetric in $x_1$, i.e., for every $x
  \in \Omega$ and $s \in
  [-1,1]$ we have $(sx_1,x_2,\dots,x_N) \in \Omega$.  
\item[(D2)] For every $\lambda>0$, the set $\Omega_\lambda:= \{x \in \Omega\::\: x_1 >\lambda\}$ has at most finitely many connected components. 
\item[(F1)] $f: (0,\infty)\times \Omega\times \cB \to \R$ is
  continuous. Moreover, for every bounded subset $K \subset \cB$ there exists $L=L(K)>0$ such
  that $\sup \limits_{x\in\Omega,\, t> 0} |f(t,x,u)-f(t,x,v)|\leq L
 |u-v|$ for $u,v \in K$.
 \item[(F2)] $f$ is symmetric in $x_1$ and nonincreasing in $|x_1|$, i.e., 
for every $t \in (0,\infty)$, $u \in \cB$, $x \in \Omega$ and $s \in
[-1,1]$ we have $f(t,s x_{1},x_2,\dots,x_N,u) \ge f(t,x,u)$.
\end{enumerate}

We note that $(D2)$ is a technical assumption which is needed for some
but not all of our results. The main result of this paper is the following.

\begin{theorem}\label{sec:goal}
 Let (D1), (F1), (F2) be satisfied, and let $u$ be a nonnegative
 global solution of $(P)$ satisfying the following conditions:
\begin{enumerate}
 \item[(U1)] There is $c_{u}>0$ such that $\|u(t)\|_{L^{\infty}}\leq c_{u}$ for every $t>0$.
\item[(U2)] The functions $u(\tau+\cdot,\cdot)$, $\tau \geq 1$ are
  uniformly equicontinuous on $[0,1]\times \overline \Omega$, that is
$$
 \lim_{h\to 0} \sup_{\substack{\tau \ge 1\\ x,\tilde{x}\in\overline{\Omega},\, t,\tilde{t}\in [\tau,\tau+1], \\ |x-\tilde{x}|,|t-\tilde{t}|<h}}|u(t,x)-u(\tilde{t},\tilde{x})|=0.
$$
\end{enumerate}
Suppose in addition that $(D2)$ holds\ \underline{or}\ that $z \not \equiv 0$ for
every $z \in \omega(u)$.\\
Then $u$ is asymptotically symmetric in $x_1$, i.e., for all $z\in \omega(u)$ we have
$z(-x_{1},x')=z(x_{1},x')$ for all $(x_{1},x')\in \Omega$.\\
Moreover, for every $z \in \omega(u)$ we have the following
alternative: Either $z \equiv 0$ on $\Omega$,
or $z$ is strictly decreasing in $|x_1|$ and therefore strictly
positive in $\Omega$.
\end{theorem}

We immediately deduce the following corollary for equilibria and
time-periodic solutions.

\begin{corollary}
\label{maincoro-simple}
 Let (D1) be satisfied for $\Omega$.
\begin{itemize}
\item[(i)] Let $f:\overline \Omega \times \cB \to
\R$, $(x,u) \mapsto f(x,u)$ be such that 
\begin{itemize}
\item[(i.1)] $f$ is continuous in $x \in \overline \Omega$ and
locally Lipschitz in $u$ uniformly with respect to $x$;
\item[(i.2)] $f$ is symmetric in $x_1$ and nonincreasing in $|x_1|$, i.e., 
for every $u \in \cB$, $x \in \Omega$ and $s \in
[-1,1]$ we have $f(s x_{1},x_2,\dots,x_N,u) \ge f(x,u)$.
\end{itemize}
Moreover, let $u \in
C_{0}(\Omega)\cap\mathcal{H}^{s}_0(\Omega)$ be a nonnegative
nontrivial weak solution of the elliptic
problem  
\begin{equation}
\label{(P-simple-elliptic)}
(-\Delta)^s u = f(x,u) \quad \text{in $\Omega$},\qquad \quad 
    u = 0 \quad \text{on $\R^N \setminus \Omega$},
\end{equation}
i.e., we have $u(x) \in \cB$ for a.e. $x \in \Omega$ and $\cE(u,\phi) = \int_\Omega f(x,u(x))\phi(x)\,dx$ for every $\phi
  \in \mathcal{H}^{s}_0(\Omega).$ Then $u$ is symmetric in $x_1$ and strictly decreasing in $|x_1|$.
\item[(ii)] Suppose that $f: (0,\infty)\times \Omega\times \cB \to \R$
  satisfies (F1), (F2) and is periodic in $t$, i.e. there is $T>0$
such that $f(t+T,x,u)=f(t,x,u)$ for all $t,x,u$. Suppose furthermore
that $u$ is a nontrivial nonnegative $T$-periodic solution of $(P)$, i.e.,
$u(t+T,x)=u(t,x)$ for all $x\in \Omega, t \in (0,\infty)$. Suppose finally
that either $(D2)$ holds or that $u(t,\cdot) \not \equiv 0$ on $\Omega$ for
all $t$. Then $u(t,\cdot)$ is 
symmetric in $x_1$ and strictly decreasing in $|x_1|$
for all times $t \in (0,\infty)$.  
\end{itemize}
\end{corollary}

\begin{remark}\label{existence}
(i) The nonnegativity assumption on $u$ in Theorem~\ref{sec:goal} can
be weakened in special cases. More precisely, if the other assumptions of
Theorem~\ref{sec:goal} are satisfied, $u(t_0,\cdot)$ is nonnegative on $\Omega$ for some $t_0>0$ and
$f(t,\cdot,0) \ge 0$ for all $t\ge t_0$, then $u(t,\cdot)$ is
nonnegative for $t \ge t_0$ as a consequence of the weak maximum principle in the form discussed in
Remark~\ref{sec:small-volume-maximum} below. Thus
Theorem~\ref{sec:goal} applies to $u$ after a time shift.\\ 
(ii) Assumption $(U2)$ implies that $\{u(t,\cdot):t>0\}\subset
C_{0}(\overline{\Omega})$ is relatively compact and therefore
$\omega(u)$ is nonempty. In Proposition~\ref{sec:appendix-6} below we
give sufficient conditions for $(U2)$ to hold.\\
(iii) In the case where, in addition to the assumptions of 
Theorem~\ref{sec:goal}, $\Omega \subset \R^N$ is a ball
centered at zero and $f$ is radially symmetric, i.e. $f(t,x,u)= \tilde
f(t,|x|,u)$, it follows -- by the invariance of the equation under
rotations -- from Theorem~\ref{sec:goal} that every $z\in \omega(u)$
is radially symmetric as well. In the special case of equilibria,
i.e., solutions of (\ref{(P-simple-elliptic)}), this has been proved in \cite{BLW05} under more restrictive assumptions on
the nonlinearity.\\
 (iv) We point out that we do not require an a priori 
 positivity assumption on elements in $\omega(u)$ in
 Theorem~\ref{sec:goal}, and thus we also do not need to assume strict
 positivity of solutions of (\ref{(P-simple-elliptic)}) in
 Corollary~\ref{maincoro-simple}. This is a special feature of the
 nonlocal problems $(P)$ and (\ref{(P-simple-elliptic)}). The strong maximum principle given by
 Theorem~\ref{sec:goal} for elements $z \in \omega(u)$ and by
 Corollary~\ref{maincoro-simple} for nonnegative solutions of
 (\ref{(P-simple-elliptic)}) is a consequence of the {\em monotonicity} 
of the nonlinearity, and it is derived as a byproduct of the method
proving the symmetry results (see in particular
Lemma~\ref{sec:schritt2} below). This contrasts with the local case $s=1$,
where counterexamples show that such a strong maximum principle is
false, see \cite[Theorem 1.1]{PT12}, \cite[Section 5]{P07} and the
references therein. In this case, an additional positivity
assumption as e.g. in  \cite[Theorem 2.2]{P07} is necessary to obtain asymptotic symmetry.
\end{remark}

The proof of Theorem~\ref{sec:goal} is
based on a parabolic variant of the moving plane method. As far as the
main structure of the argument is concerned, we follow the strategy elaborated
by Pol\'a\v cik \cite{P07,P09} in the context of Dirichlet problems for fully nonlinear parabolic
differential equations, but we need new and quite different tools. We recall that the moving plane method has its
roots in a classical work of Alexandrov \cite{A62} on constant mean
curvature surfaces and Serrin \cite{serrin} on overdetermined boundary
value problems, whereas Gidas, Ni and Nirenberg \cite{GNN79} provided
the framework to consider Dirichlet problems
for nonlinear elliptic differential equations. In the case where the
underlying domain is $\R^N$, the method of moving plane has been
applied in integral form in \cite{Chen_Li_Ou,felmer-quaas-tan}
to deduce symmetry and classification results for solutions of
semilinear elliptic equations involving the fractional Laplacian.  
Birkner, Lop\'ez-Mimbela and Wakolbinger \cite{BLW05} used a variant of the
moving plane method, paired with probabilistic methods, to prove
radial symmetry of all equilibria of $(P)$ in the case where the
underlying domain is the unit ball $B$ and the nonlinearity $f$ is
nonnegative, independent of $t$ and $x$, and nondecreasing in $u$. Up to the authors' knowledge, our results are the first symmetry
results for parabolic boundary value problems involving the fractional
Laplacian and even for the elliptic problem 
if $f$ depends on $x$ or the domain is more general than a ball. We
point out that -- in comparison with the elliptic case -- proving asymptotic
symmetry in the parabolic setting with the moving plane approach 
requires much finer -- time dependent -- estimates. This is already
evident from the seminal work of Pol\'a\v
cik \cite{P09,P07} for the case of nonlinear differential equations.
One key requirement is a special version of a parabolic Harnack
inequality related to a linear fractional diffusion equation. Felsinger and Kassmann derived a parabolic Harnack inequality in \cite{FK12}, which
requires nonnegativity of the solutions in the entire space. This
global nonnegativity assumption is not technical since -- already in
the elliptic case -- 
the Harnack inequality for the fractional Laplacian is not valid in a
purely local form, see e.g. \cite[Theorem 2.2]{K11} for a counterexample.
However, since the moving plane method consists in studying the difference
between the reflection of a solution of $(P)$ at a hyperplane and the
solution itself, we need to derive a corresponding Harnack
inequality for antisymmetric (and therefore sign changing) supersolutions of a
class of linear problems in the present paper. Another (closely related) problem in the fractional setting is the lack of
local comparison principle to derive estimates via sub- or
supersolutions. Here much finer quantitative arguments are needed
to control the nonlocal effects and exclude the appearance of
intersections in finite time. We will establish such estimates in two
steps in Section~\ref{sub} below, passing first to the
Caffarelli-Silvestre extension of the solution $u$, which is defined,
for each fixed time, on the half space $\R^{N+1}_+$ (see \cite{CS07}). 

It seems worthwhile to note that another type of Dirichlet boundary
conditions has also been assigned to the fractional Laplacian in
the literature. In \cite{Tan,CDDS11,CT10}), the authors consider
the $s$-th power of the Dirichlet Laplacian in spectral theoretic
sense, which -- in the case of a bounded domain $\Omega$ -- 
is given by $A^{s}_{\Omega} u := \sum \limits_{k=1}^{\infty}\mu_{k}^{s}u_{k}e_{k}$.
Here $\mu_{k}=\mu_{k}(\Omega)$ are the eigenvalues of the Dirichlet
Laplacian on $\Omega$ in increasing order (counted with multiplicity),
$e_{k}$, $k \in \N$ are the corresponding eigenfunctions and
$u_{k}:=\int_{\Omega} u e_k\ dx$ the corresponding Fourier
coefficients of $u$. 

In order to explain the role of $A^{s}_{\Omega}$ in terms of stochastic processes, we recall that the
$2s$-stable process is constructed by subordinating Brownian
motion with a $s$-stable subordinator, see \cite[Chapter 1.3]{A09}. 
On the other hand, the process generated by $A^{s}_{\Omega}$ is obtained by first killing Brownian
motion upon leaving $\Omega$ and then subordinating this process with
a $s$-stable subordinator, see e.g. \cite{SV03}. Hence the order of killing and
subordination is reversed in this case. It is easy to see that the
corresponding operators coincide only if $\Omega=\R^N$ (where the
Dirichlet boundary conditions are not present). For more information related
to these stochastic processes and their generators, we refer the
reader to \cite{BKK08},\cite{J05} or \cite[Chapter 3]{A09}.
It is natural to ask whether a result similar to Theorem~\ref{sec:goal} is true for the
corresponding problem with the operator $A^s_{\Omega}$. For elliptic semilinear problems
involving the operator $A^s_{\Omega}$, symmetry and monotonicity results 
 have been proved recently in special cases in
\cite{CT10,CDDS11} by applying the moving plane method to the 
Caffarelli-Silvestre extensions of the solutions.  

The article is organized as follows. In Section~\ref{1}, we develop
the new tools we need to carry out the moving
plane method for the fractional parabolic problem $(P)$. We believe
that the results of this Section could be of interest for other
problems as well. Since, as already noted, the moving plane method consists in studying the difference
between the reflection of a solution of $(P)$ at a hyperplane and the
solution itself, we are led to study antisymmetric supersolutions of linear problems here. Due to the nonlocality of the 
fractional Laplacian, it is important to
estimate the influence of the negative part of these functions. This is one of the key differences
in comparison with local problems involving classical differential
operators. The first part of this Section is concerned with a
parabolic small volume maximum principle. In Section~\ref{stochastik}
we establish, based on recent results in \cite{FK12}, a parabolic Harnack inequality for antisymmetric
supersolutions of a class of linear fractional problems. Section~\ref{sub} is devoted
to a generalized subsolution estimate. The idea to control the positive part of the
solution by comparing with suitable subsolutions is inspired by
\cite{P07}. However, as mentioned above, the argument is
essentially more involved in the present setting, and this is the only
stage where we had to pass to the
Caffarelli-Silvestre extension. In
Section~\ref{mainlinear}, we combine all estimates obtained so far
to deduce our main result on antisymmetric supersolutions
for a class of linear problems. This result should be seen as an analogue of
\cite[Theorem 3.7]{P07} for the fractional case. The moving plane argument is
then carried out in Section~\ref{sec:proof-main-symmetry}. Here we
follow the main structure of the argument in \cite[Chapter 4]{P07},
but we need to implement some new ideas at key points (see in
particular the proof of Lemma~\ref{sec:schritt2}) in the nonlocal
setting. In the appendix, we present a sufficient condition for 
$(U2)$, and we discuss a specific example to which
Theorem~\ref{sec:goal} applies.

\subsection{Notation}
\label{sec:notation}

The following notation is used throughout the paper. For $x \in \R^N$
and $r>0$, $B_r(x)$ is the open ball centered at $x$ with radius
$r$ and $\omega_{N}$ will denote the volume of the $N$-dimensional ball with radius $1$.
For any subset $M \subset \R^N$, we denote by $1_M: \R^N \to \R$ the
indicator function of $M$ and $\diam(M)$ the diameter of
$M$. Moreover, we let $\inrad(M)$ denote the supremum of
all $r>0$ such that \textit{every} %every kursiv um die Änderung zu verdeutlichen 16.11.
connected component of $M$ contains a ball 
$B_{r}(x_{0})$ with $x_0 \in M$. This notation -- taken from
\cite{P07} -- differs slightly
from the usual one but is very convenient in our setting. If $T \subset \R$, $\Omega \subset
\R^N$ are subsets and $u: T \times \Omega \to \R$, $(t,x) \mapsto u(t,x)$
is a function, we frequently write $u(t)$ in place of $u(t,\cdot): \Omega
\to \R$ for $t \in T$. If $M \subset \R^{N}$ resp. $M \subset
\R^{N+1}$ is a subset and $w: M \to \R$ is a function, the
inequalities $w \ge 0$ and $w>0$ are always understood in pointwise
sense. Moreover, $w^+= \max\{w,0\}$ resp. $w^-=-\min\{w,0\}$ denote the
positive and negative part of $w$, respectively. If $M$ is measurable with $|M|>0$ (where $|\cdot|$
always stands for Lebesgue measure) and $w \in L^1(M)$, we put 
$$
[w]_{L^1(M)}:= \frac{1}{|M|}\int_{M}w(x)\,dx,\quad   [w]_{L^1(M)}:=
\frac{1}{|M|}\int_{M}w(t,x)\,dt dx ,\quad \text{respectively,}
$$
to denote the mean of $w$ over $M$. If $D,U \subset \R^N$ are subsets,
the notation $D \subset \subset U$ means that $\overline D$ is compact
and contained in the interior of $U$. Moreover, we set 
$$
\dist(D,U):= \inf\left\{|x-y|\;:\; x \in D,\, y \in U\right\},
$$
so this notation does {\em not} stand for the usual Hausdorff
distance. If $D= \{x\}$ is a singleton, we simply write $\dist(x,U)$
in place of $\dist(\{x\},U)$. Finally, when we call an interval $T \subset \R$ a time interval, we
assume that it consists of more than one point.\\[0.1cm]

\section{Antisymmetric supersolutions of a corresponding linear problem}\label{1}
\label{sec:antisymm-supers-line}

Throughout this section, we consider a fixed open half space $H$ and the
reflection $Q: \R^N \to \R^N$ at $\partial H$. We will
call a function $w: \R^N \to \R$ antisymmetric if $w(Q(x))=-w(x)$ for
every $x \in \R^N$, i.e., $w$ is antisymmetric with respect to $Q$.
We first fix notions of supersolutions. For an open subset $U' \subset
\R^N$, we introduce the function space 
\begin{equation}
  \label{eq:66}
\cV^s(U'):= \{u \in L^\infty(\R^N)\::\: \int_{U' \times U'}
\frac{|u(x)-u(y)|^2}{|x-y|^{N+2s}}\,dxdy <\infty\},
\end{equation}
endowed with the norm 
$$
\|u\|_{\cV^s(U')}:= \|u\|_{L^\infty(\R^N)}+ \Bigl(\int_{U' \times U'}
\frac{|u(x)-u(y)|^2}{|x-y|^{N+2s}}\,dxdy\Bigr)^{\frac{1}{2}}.
$$
We note that if $U \subset \subset U'$ is a pair of open sets and $u
\in \cV^s(U')$, $v \in \mathcal{H}^s_0(U)$, then 
$\cE(u,v)$ is well defined by (\ref{eq:46}).

\begin{definition}
\label{sec:antisymm-supers-corr-2}
Let $U \subset \R^N$ be a bounded open subset, $T$ a time interval and $c,g \in
L^\infty(T \times U)$. We call a function $v: T \times \R^N \to \R$ a
{\em supersolution} of 
\begin{equation}
  \label{eq:3-g}
\partial_t v +\hl v = c(t,x) v +g(t,x)
\end{equation}
on $T \times U$ if $v \in C(T,\cV^{s}(U'))\cap
C^{1}(T,L^{2}(U))$ for some open set $U' \subset \R^N$ with $U \subset
\subset U'$ and 
$$
\mathcal{E}(v(t),\phi) \ge \int_{U} (c(t,x)v(t)+g(t,x)-\partial_{t}
v(t))\varphi\ dx
$$
for all $\phi \in
\mathcal{H}^{s}_{0}(U)$, $\phi \ge 0$ and a.e. $t\in T$. If, in addition, $U \subset H$ and $v$ is antisymmetric, we call $v$ an {\em antisymmetric
supersolution}. A supersolution of (\ref{eq:3-g}) on $T \times U$ will be called an {\em entire supersolution} if $v \ge
0$ on $T \times (\R^N \setminus U)$. If $U \subset H$, an
antisymmetric supersolution of (\ref{eq:3-g}) on $T \times U$ will be
called an {\em entire antisymmetric supersolution} if $v \ge
0$ on $T \times (H \setminus U)$. 
\end{definition}

\begin{remark}
\label{sec:antisymm-supers-corr-4}  
(i) Note that an entire antisymmetric supersolution
$v$ of (\ref{eq:3-g}) on $T \times U$ may take negative values in $\R^N
\setminus H$, so in general it is {\em not} an entire
supersolution of (\ref{eq:3-g}).\\
(ii) Let $T,U$ and $c$ be as in the definition above. We will mostly
consider the case $g \equiv 0$ in the remainder of the paper, i.e., we
consider supersolutions of 
\begin{equation}
  \label{eq:3-0}
\partial_t v +\hl v = c(t,x) v
\end{equation}
on $T \times U$. We briefly explain the connection between $(P)$ and
(\ref{eq:3-0}). Suppose that $(F1)$ is satisfied and that 
\begin{align}
&H \cap \Omega \not= \varnothing, \qquad Q(H \cap \Omega) \subset
\Omega \qquad \text{and} \label{eq:39}\\ 
&f(t,Q(x),u) \ge f(t,x,u)\qquad \text{for every $t \in (0,\infty),x \in
  U$ and $u \in \cB$.}  \label{eq:22}
\end{align}
Let $u$ be a nonnegative solution
of $(P)$, and let $v(t,x)= u(t,Q(x)) -
u(t,x)$ for $x \in \R^N,\,t \ge 0$. Then $v$ is an entire antisymmetric supersolution of (\ref{eq:3-0}) with
$T=(0,\infty)$, $U= H \cap \Omega$ and  
$$
c(t,x)= \left \{
  \begin{aligned}
  &\frac{f(t,x,u(t,Q(x))-f(t,x,u(t,x))}{v(x)},&&\qquad u(t,Q(x)) \not= u(t,x);\\
  &0, &&\qquad u(t,Q(x))= u(t,x). 
  \end{aligned}
\right.
$$
Indeed, by (\ref{eq:39}) we have $v \ge 0$ on $T \times (H \setminus
U)$. Moreover, for $\phi \in \mathcal{H}^{s}_{0}(U)$, $\phi
\ge 0$ and $t \in (0,\infty)$ we have 
\begin{align*}
\cE(v(t),\phi)&= \cE(u(t) \circ Q - u(t),\phi)=\cE(u(t),\phi \circ Q - \phi) \\
&= \int_{\Omega}(f(t,x,u)-\partial_{t}
u)[\varphi \circ Q - \varphi]\,dx \\
&=\int_{U}[f(t,Q(x),u(t,Q(x)))-f(t,x,u(t,x)) -\partial_{t}
(u \circ Q - u)]\varphi\,dx\\
& \ge \int_{U}[c(t,x) v - \partial_{t}v]\varphi\,dx,
\end{align*}
where (\ref{eq:22}) was used in the last step. 
\end{remark}

The following
observation will be useful in the sequel. 

\begin{lemma}
\label{sec:antisymm-supers-corr-1}
For any $\phi \in \cH_0^s(H)$ and every antisymmetric $v \in
H^s(\R^N)$ we have 
\begin{equation}
  \label{eq:37}
\mathcal{E}(v,\varphi)= \frac{1}{2} \int_{H}\int_{H}(v(x)-v(y))({\varphi}(x)-{\varphi}(y))J(x,y) dx\ dy+ 2\int_{H} \kappa_{H}(x) v(x)\varphi(x)\ dx 
\end{equation}
with 
\begin{align*}
J(x,y)&=  \frac{c_{N,s}}{|x-y|^{N+2s}}
  -\frac{c_{N,s}}{|x-Q(y)|^{N+2s}}\\
\text{ and } \kappa_H(x)&= \int_{\R^N \setminus H}\frac{c_{N,s}}{|x-y|^{N+2s}}\,dy =
\frac{ 4^s \Gamma(\frac{1}{2}+s) }{\sqrt{\pi} \Gamma(1-s)} [\dist(x,\partial H)]^{-2s}
\end{align*}
for $x,y \in H$, where $c_{N,s}$ is given in (\ref{eq:45}). Moreover,
\begin{equation}
  \label{eq:38}
J(x,y) \geq\frac{c_{N,s}[1-5^{-N/2-s}]}{|x-y|^{N+2s}}
\end{equation}
for $x,y \in H$ with $|x-y| \le \min \{\dist(x,\partial H),
\dist(y,\partial H)\}$.
\end{lemma}

\begin{proof}
It is convenient to write $\bar x$ in place of $Q(x)$ for $x \in \R^N$
in the following. For $\phi \in \cH_0^s(H)$ and an antisymmetric $v \in
H^s(\R^N)$ we then have 
\begin{align*}
\mathcal{E}(v,\varphi) &=
\frac{c_{N,s}}{2}\biggl(\int_{H}\int_{H}\frac{(v(x)-v(y))({\varphi}(x)-{\varphi}(y))}{|x-y|^{N+2s}}
  \ dx\ dy\\
&\qquad\qquad\qquad\qquad\qquad\qquad + \int_{H}\int_{\R^N \setminus H}\dots dx\ dy + \int_{\R^N \setminus H}\int_{H}\dots dx\ dy\biggr) \\
&= \frac{c_{N,s}}{2}\int_{H}\int_{H}\Bigl[
\frac{(v(x)-v(y))({\varphi}(x)-{\varphi}(y))}{|x-y|^{N+2s}} \\
&\qquad\qquad\qquad\qquad\qquad\qquad-
\frac{(v(\bar x)-v(y))\varphi(y)}{|\bar x-y|^{N+2s}}+\frac{(v(x)-v(\bar
y))\varphi(x)}{|x-\bar y|^{N+2s}}\Bigr]\ dx \ dy\\
&=
\frac{c_{N,s}}{2}\int_{H}\int_{H}\frac{(v(x)-v(y))({\varphi}(x)-{\varphi}(y))}{|x-y|^{N+2s}}
  \ dx\ dy  \\
&\qquad\qquad\qquad\qquad\qquad\qquad+ c_{N,s}
  \int_{H} \int_{H}\frac{(v(x)+v(y))\varphi(y)}{|x-\bar y|^{N+2s}}\ dx \ dy\\
&=\frac{1}{2}\int_{H}\int_{H}(v(x)-v(y))({\varphi}(x)-{\varphi}(y))J(x,y)\
dx\ dy \\
&\qquad\qquad\qquad\qquad\qquad\qquad+ 2c_{N,s} \int_{H}\int_{H}\frac{v(y)\varphi(y)}{|x-\bar y|^{N+2s}}\ dx\ dy\\
&=\frac{1}{2} \int_{H}\int_{H}(v(x)-v(y))({\varphi}(x)-{\varphi}(y))J(x,y) dx\ dy+ 2\int_{H} \kappa_{H}(x) v(x)\varphi(x)\ dx 
\end{align*}
with $J$ and $\kappa_H$ as defined above, as claimed. 
% Where the constant can be simply calculated via
% \begin{align*}
%  \int_{\R^N \setminus H}\frac{c_{N,s}}{|x-y|^{N+2s}}\,dy&=c_{N,s}\int_{x_{1}}^{\infty}\int_{\R^{N-1}}(y_{1}^{2}+|y'|^2)\frac{-N/s-s}\,dy'dy_1\\
% &=c_{N,s}(N-1)\omega_{N-1}\int_{x_{1}}^{\infty}\int_{0}^{\infty}(y_{1}^{2}+r^2)\frac{-N/s-s}r^{N-2}\,drdy_1\\
% &=\frac{c_{N,s}(N-1)\pi^{\frac{N-1}{2}}\Gamma(1/2+s)\Gamma(\frac{N-1}{2})}{2\Gamma(\frac{N+1}{2})\Gamma(N/2+s)}2\int_{x_{1}}^{\infty}y_{1}^{-2s-1}\,dy_1\\
% &=\frac{c_{N,s}\pi^{\frac{N-1}{2}}\Gamma(1/2+s)}{\Gamma(N/2+s)s}x_1^{-2s}=\frac{4^s\Gamma(1/2+s)}{\sqrt{\pi}\Gamma(1-s)}\dist(x,\partial H)^{-2s}
% \end{align*}
To see
(\ref{eq:38}), let
$d>0$ and $x,y
  \in H$ with $|x-y| \le d \le \min \{\dist(x,\partial H),
  \dist(y,\partial H)\}$. Then $|x-\bar{y}|^2\geq
|x-y|^2+4d^2$ and therefore 
$$
\frac{|x-y|^2}{|x-\bar{y}|^2} \leq \frac{|x-y|^2}{|x-y|^2+4d^2} \leq \frac{1}{5},
$$
which implies that 
$$
\frac{J(x,y)|x-y|^{N+2s}}{c_{N,s}} =
\left(1-\left(\frac{|x-y|^2}{|x-\bar{y}|^2}\right)^{\frac{N+2s}{2}}\right)
\geq 1-5^{-N/2-s}
$$
as claimed in (\ref{eq:38}).
\end{proof}

\subsection{A small volume maximum principle}
\label{sec:antisymm-supers-corr}
The main result of this subsection is the following.

\begin{proposition}
\label{sec:klein1-1}
For every $c_\infty, \gamma>0$ there exists
$\delta=\delta(N,s,c_\infty,\gamma)>0$ such that for any bounded open subset $U
\subset H$ with $|U| \le \delta$, any time interval $T:=[t_{0},t_1]$, any $c \in L^\infty(T \times U)$ with 
$\|c^+\|_{L^\infty} \le c_\infty$ and any entire antisymmetric supersolution $v$ 
of (\ref{eq:3-0}) on $T \times U$ we have 
\begin{equation}\label{sec:kleinklein1}
 \|v^{-}(t,\cdot)\|_{L^{\infty}(H)} \leq 
 e^{-\gamma(t-t_{0})} \|v^{-}(t_{0},\cdot)\|_{L^\infty(H)} \qquad \text{ for all } t\in T.
\end{equation}
\end{proposition}

In the proof we will use the following standard estimate, see
e.g. \cite[Lemma 6.1]{NPV11}.

\begin{lemma}\label{mengen}
For every measurable $A\subset \R^{N}$ and every $x\in \R^{N}$ we have
\[
 \int_{\R^{N}\setminus A}\frac{1}{|x-y|^{N+2s}}\ dy\geq K |A|^{-\frac{2s}{N}}
\]
with $K=K(N,s)=\frac{N}{2s}\omega_{N}^{1+2s/N}$.
\end{lemma}

\begin{proof}[Proof of Proposition~\ref{sec:klein1-1}]
For given $c_\infty, \gamma>0$ we put $\delta:=
\Bigl(\frac{c_{N,s}\: K}{\gamma+ c_\infty}\Bigr)^{\frac{N}{2s}}$, where
$K$ is given in Lemma~\ref{mengen}. By assumption and
Lemma~\ref{mengen}, we then have 
\begin{equation}
  \label{eq:2}
\kappa_{U}(x):= \int_{\R^{N}\setminus U}\frac{c_{N,s}}{|x-y|^{N+2s}}\ dy \ge \gamma
+c_\infty \qquad \text{for every $x \in \R^N$.}
\end{equation}
Without loss of generality, we may assume that $t_{0}=0$. Let
$d:=\|v^{-}(0)\|_{L^{\infty}(U)}$, and define $u(t,x):=e^{\gamma
  t}v(t,x)$ for $t \in [0,t_1]$, $x \in \R^N$. Then $u$ is an antisymmetric
supersolution of $u_t +(-\Delta)^s u = \tilde c(t,x) u$ on $U$ with
$\tilde{c}(t,x)=c(t,x)+\gamma$. We need to show that 
\begin{equation}
  \label{eq:34}
u(t,x) \ge -d \qquad \text{for $x \in H$ and $t \in [0,t_1]$.}  
\end{equation}
For $0 \le t \le t_1$, we consider the function
$\phi(t)=\phi(t,\cdot): \R^N \to \R$ defined by
$\varphi(t,x)=\left(u(t,x)+d\right)^{-}1_{H}(x)$. Since $u(t) \in
\cV^s(U')$ for some open set $U'$ with $U \subset \subset U'$ and $u
\ge 0$ in $H \setminus U$, it follows from \cite[Lemma 5.1]{NPV11}
that $\phi(t) \in
\mathcal{H}^{s}_{0}(U)$ for $0 \le t \le t_1$. We then have
\begin{equation}\label{prob2}
 \mathcal{E}(u(t),\varphi(t)) \geq \int_{U} (\tilde c(t,x)u-\partial_{t}
u)\varphi\ dx \ge \int_{U}(c_\infty +\gamma)u(t) \varphi(t) \ dx +
\frac{1}{2}\frac{d}{dt}\int_{U}\varphi(t)^{2}\ dx. 
\end{equation}
We first claim that 
\begin{equation}
  \label{eq:33}
 \mathcal{E}(u(t),\varphi(t))\leq
- \mathcal{E}(u^{-}(t)1_{H},\varphi(t))\qquad \text{for $(t,x) \in [0,t_1]$.}
\end{equation}
Indeed, for $(t,x) \in [0,t_1] \times \R^N$ we have
\begin{align*}
\left(u(t,x)-u(t,y)\right)\left(\varphi(t,x)-\varphi(t,y)\right)+\left(u^{-}(t,x)1_{H}(x)-u^{-}(t,y)1_{H}(y)\right)\left(\varphi(t,x)-\varphi(t,y)\right)& \\
=-\Bigl[\varphi(t,x)\left(u(t,y)+u^{-}(t,y)1_{H}(y)\right)+\varphi(t,y)\left(u(t,x)
  + u^{-}(t,x)1_{H}(x)\right)\Bigr]&.
\end{align*}
Thus we find, using the symmetry of the kernel and the antisymmetry of
$u$,
\begin{align*}
\mathcal{E}&(u(t),\varphi(t))+\mathcal{E}(u^{-}(t)1_{H},\varphi(t))=-c_{N,s}\int_{\R^{N}}\varphi(t,y)
\int_{\R^{N}}\frac{\left(u(t,x) +1_{H}(x)u^{-}(t,x)\right)}{|x-y|^{N+2s}}\ dxdy \\
 &=-c_{N,s}\int_{H}\varphi(t,y)\int_{H}\Bigl(
 \frac{u^+(t,x)}{|x-y|^{N+2s}} -
 \frac{u(t,x)}{|Q(x)-y|^{N+2s}}\Bigr)dx dy
\end{align*}
and hence $\mathcal{E}(u(t),\varphi(t))+\mathcal{E}(u^{-}(t)1_{H},\varphi(t))\le
0$ for $t \in [0,t_1]$, since $|x-y| \le
|Q(x)-y|$ for $x,y \in H$ and $\phi$ is nonnegative. This shows (\ref{eq:33}). We now put 
$$
A_1(t):= \{x \in H\::\: u(t,x)
\le -d\}\qquad \text{and}\qquad A_2(t):= (\R^N \setminus H) \cup \{y \in H\::\: u(t,y)
> -d\}
$$
for $t \in [0,t_1]$. Then we have, for $t \in [0,t_1]$,
\begin{align}
\mathcal{E}(&u^{-}(t)1_{H},\varphi(t))\nonumber\\
&=\frac{c_{N,s}}{2}\int_{A_1(t)}\int_{A_1(t)}
\frac{\left(u^{-}(t,x)1_{H}(x)-u^{-}(t,y)1_{H}(y)\right)\left(\varphi(t,x)
    -\varphi(t,y)\right)}{|x-y|^{N+2s}}\,dxdy \nonumber\\
&\quad \qquad  +
c_{N,s}\int_{A_1(t)}\varphi(t,x)\int_{A_2(t)}\frac{u^{-}(t,x)1_{H}(x)-u^{-}(t,y)1_{H}(y)}{|x-y|^{N+2s}}\,dydx \nonumber\\
&\ge \frac{c_{N,s}}{2} \int_{A_1(t)}\int_{A_1(t)}  \frac{\left(\varphi(t,x)
    -\varphi(t,y)\right)^{2}}{|x-y|^{N+2s}}dxdy \nonumber\\
&\quad \qquad + c_{N,s}\int_{\R^N}\varphi(t,x)
\int_{A_2(t)}\frac{d-u^{-}(t,y)1_{H}(y)}{|x-y|^{N+2s}}\ dydx
\nonumber\\
&\ge c_{N,s}\int_{\R^N}\varphi(t,x)
\int_{\R^N \setminus U}\frac{d}{|x-y|^{N+2s}}\ dydx  \ge 
d \int_{\R^N}\varphi(t,x) \kappa_{U}(x)\, dx \label{eq:35}
\end{align}
with $\kappa_{U}(x)$ as defined in (\ref{eq:2}). Let $h(t):= \|\varphi(t,\cdot)\|_{L^{2}(U)}^{2}$ for $t \in
[0,t_1]$. Combining (\ref{prob2}),~(\ref{eq:33}) and (\ref{eq:35}), we get 
\begin{align*}
h'(t)&\leq -2(\gamma+c_\infty)\int_{U}u(t,x)\varphi(t,x)\,dx -
 2d \int_{\R^N} \varphi(t,x) \kappa_{U}(x)\, dx\\
&\le 2(\gamma+c_\infty)h(t) + 2d \int_{U}[\gamma+c_\infty-
\kappa_U(x)] \varphi(t,x)\, dx \qquad \text{for $t \in [0,t_1]$.}
\end{align*}
By (\ref{eq:2}) we conclude that $h'(t) \le 2 (\gamma+c_\infty)h(t)$
for $t \in [0,t_1]$. Since $h(0)=0$, we infer $h(t)= 0$ for $t \in [0,t_1]$. This
shows (\ref{eq:34}), as required.
\qedhere
\end{proof}

\begin{remark}
\label{sec:small-volume-maximum}
For entire supersolutions $v$ of (\ref{eq:3-0}), a
corresponding small volume maximum principle can be derived in a
similar but much easier way. More precisely, for every $c_\infty, \gamma>0$ there exists
$\delta>0$ such that for any bounded open subset $U
\subset \R^N$ with $|U| \le \delta$, any time interval $T:=[t_{0},t_1]$, any $c \in L^\infty(T \times U)$ with 
$\|c^+\|_{L^\infty} \le c_\infty$ and any entire supersolution $v$ 
of (\ref{eq:3-0}) on $T \times U$ we have 
\begin{equation}\label{sec:kleinklein1-1}
 \|v^{-}(t,\cdot)\|_{L^{\infty}(U)} \leq 
 e^{-\gamma(t-t_{0})} \|v^{-}(t_{0},\cdot)\|_{L^\infty(U)} \qquad \text{ for all } t\in T.
\end{equation}
As a consequence, we may readily derive the following weak maximum principle: If $\Omega
\subset \R^N$ is a bounded open subset, $T:=[t_{0},t_1]$ a time interval, $c \in L^\infty(T \times \Omega)$ and $v$ an entire supersolution of
(\ref{eq:3-0}) on $T \times \Omega$ such that $v(t_0,x) \ge 0$ for
a.e. $x \in \Omega$, then also $v(t,x)\geq 0$ for all $t \in T$ and
almost every $x \in \Omega$. 
\end{remark}

\subsection{A Harnack inequality for antisymmetric
      supersolutions}\label{stochastik}

In this part we state a Harnack inequality for antisymmetric
supersolutions of (\ref{eq:3-0}). We will derive this inequality -- via a
reformulation of the problem -- from a
recent result in \cite{FK12}. We
need to introduce some notation.
Denote by $\vartriangle = \{(x,x)\::\: x \in \R^N\}$ the diagonal in
 $\R^N \times \R^N$. We fix $r_0 \in (0,1]$ and $C_1,C_2>0$, and
we consider a function $k: \R^N \times
\R^N \setminus \vartriangle \to [0,\infty)$ satisfying, for every $x,y \in \R^N$ with $x \not=y$,
\begin{equation}
  \label{eq:14}
  \begin{aligned}
&k(x,y)=k(y,x);\\ 
&k(x,y) \leq C_{1}|x-y|^{-N-2s};\\
&k(x,y) \ge C_2|x-y|^{-N-2s}\quad \text{if $|x-y|\leq r_0$.}    
  \end{aligned}
\qquad 
\end{equation}
The quadratic form corresponding to this kernel is given by 
$$
\cE_k(u,v)= \frac{1}{2} \int_{\R^N} \int_{\R^N}(u(x)-u(y))(v(x)-v(y))k(x,y)\,dx dy, \text{ for $u,v \in H^s(\R^N)$.}
$$
Recall the definition of $\cV^s(U')$ in (\ref{eq:66}). If $U \subset
\R^N$ is a bounded open subset, $T \subset \R$ a time interval
with nonempty interior and $g \in L^\infty(T \times U)$, we say that a
function $v$ is a supersolution of the problem 
\begin{equation}
  \label{eq:10}
 \partial_{t}v(t,x)-P.V.\int_{\R^{N}}(v(t,y)-v(t,x))k(x,y)\ dy = g 
\end{equation}
on $T \times U$ if $v \in C(T,\cV^s(U')) \cap C^1(T,L^2(U))$ for some
open set $U' \subset
\R^N$ with $U \subset \subset U'$ and  
$$
\cE_k(v(t),\phi) \ge \int_{U}[g(t,x)- \partial_t v(t,x)]
\phi(x)\,dx
$$
for $\phi \in
  \H^s_0(U)$, $\phi \ge 0$ and a.e. $t \in T$. Next we introduce notation for parabolic cylinders. For
$t_0 \in \R$, $x_0 \in \R^N$, $r,\vartheta>0$ we put
$Q(r,\vartheta,t_{0},x_{0}):=(t_{0},t_0+8 \vartheta)\times
B_{2r}(x_{0})$ and 
$$
Q^{-}(r,\vartheta,t_{0},x_{0}):=(t_{0},t_0+\vartheta)\times
B_{r}(x_{0}),\quad Q^{+}(r,\vartheta,t_{0},x_{0}):=(t_{0}+7
\vartheta,t_{0}+8 \vartheta)\times B_{r}(x_{0}).
$$
In view of the scaling properties of (\ref{eq:10}), the following is a mere
reformulation of a special case of \cite[Theorem
1.1]{FK12}, see also \cite[Remark after Theorem 1.2]{FK12}. We point out that the notion of supersolution considered
in \cite{FK12} is weaker than the one considered here.

\begin{theorem}\label{t:harnack}
Let $r_0 \in (0,1]$ and $\vartheta,C_1,C_2>0$ be given. Then
there are constants $c_i>0$, $i=1,2$ depending on $N,s,r_0,\vartheta, C_1,C_2$ such that for any $k: \R^N \times
\R^N \setminus \vartriangle \to [0,\infty)$ satisfying (\ref{eq:14}),
any $(t_0,x_0) \in \R^{N+1}$, any $g \in L^\infty(Q(r_0,\vartheta,t_{0},x_{0}))$ and any supersolution $v$ of
(\ref{eq:10}) on $Q(r_0,\vartheta,t_{0},x_{0})$ which is nonnegative in
$(t_{0},t_{0}+8\vartheta) \times
\R^{N}$ we have
\begin{equation}
\label{eq:15}
\inf_{(t,x)\in Q^{+}(r_0,\vartheta, t_{0},x_{0})} v \geq
c_1[v]_{L^{1}(Q^{-}(r_0,\vartheta,t_{0},x_{0}))}-c_2
\|g\|_{L^\infty(Q(r_0,\vartheta,t_{0},x_{0}))}.
\end{equation}
\end{theorem}

By an argument based on building chains of cylinders, we deduce the following
Harnack inequality for general pairs of domains. We include the proof here since the argument is not
completely standard. A similar argument has been detailed in
\cite[Appendix]{P07}, but we need to argue somewhat differently since the
triples of parabolic cylinders in Theorem \ref{t:harnack} have a
smaller overlap than the ones considered in \cite{P07}.  

\begin{corollary}\label{harnack-gen}
 Let $r_0 \in (0,1]$, $R,\tau,\eps>0$ and $C_1,C_2>0$ be given. Then there exist
 positive constants \mbox{$c_i=c_i(N,s,r_0,C_1,C_2,R,\eps,\tau)>0$},
 $i=1,2$ with the following property:\\
Let $k: \R^N \times \R^N \setminus \vartriangle \to [0,\infty)$ satisfy 
(\ref{eq:14}), and let $D\subset\subset U\subset \R^N$ be a pair of
bounded domains
such that $\dist(D,\partial
 U) \geq 2 r_0$, $|D| \ge \eps$ and $\diam(D)\le R$. Moreover, let
$g \in
 L^\infty(T \times U)$ and a supersolution $v$ of
(\ref{eq:10}) on $T \times U$ be given such that $v$ is nonnegative in
$T \times \R^N$, where $T=[t_0,t_0 + 4\tau]$ for some $t_0 \in \R$. Then we have
 \begin{equation}
\inf_{(t,x)\in T_+\times D}v(t,x)\geq c_1 [v]_{L^{1}(T_- \times
  D)}-c_2 \|g\|_{L^\infty(T \times U)},
\end{equation}
where $T_+=[t_0+3 \tau,t_0 + 4 \tau]$ and
$T_-=[t_0+\tau,t_0+2 \tau]$.
\end{corollary}

\begin{proof}
We first note that there exist $n=n(N,R,r_0) \in \N$ and
$\mu=\mu(N,R,r_0)>0$ such that the following holds:\\ 
{\em For every subset $D \subset \R^N$ with $\diam\, D \le R$ there
  exists a subset $S_D \subset D$ of $n+1$ points such that $D$ is
  covered by the balls $B_{r_0}(x)$, $x \in S_D$, and for every two
points $x_*,x^* \in S_D$ there exists a finite sequence $x_j \in S_D$, $j=0,\dots,n$ such that} 
\begin{equation}
  \label{eq:29}
x_{0}=x_*,\quad x_{n}=x^* \qquad \text{and} \qquad |B_{r_0}(x_j)
  \cap B_{r_0}(x_{j+1})| \ge \mu \quad \text{for
      $j=0,\dots,n-1$.}
\end{equation}
We now fix $D \subset \subset U \subset H$ as in the assertion, and we
fix
$n,\mu$ and a set $S_D$ with the property above. Next, we
put $\vartheta= \frac{\tau}{7} \min \{\frac{1}{17},\frac{1}{n+3}\}$, and we claim the
following:\\ 
{\em For given $t_* \in [t_0+\tau,t_0+2\tau]$ and $t^* \in [t_0+3 \tau,t_0 +4
  \tau]$ there exists a finite sequence $t_*=s_0<...<s_m=t^*-8 \vartheta$ such that} 
\begin{equation}
  \label{eq:27}
s_j + 7\vartheta \le s_{j+1} \le s_j + \frac{15}{2} 
\vartheta \qquad \text{for $j=0,\dots,m-1$}
\end{equation}
and
\begin{equation}
  \label{eq:28}
\max\{14,n\} \le m \le \max\{51,3(n+3)\} 
\end{equation}
Indeed, let $m \in \N$ and $\sigma \in
[0,7 \vartheta)$ be such that $t_* + 7 m \vartheta +\sigma=t^*-
8 \vartheta$. The definition of $\vartheta$ and the restrictions on
$t_*$, $t^*$ then force (\ref{eq:28}), and (\ref{eq:27}) holds with $s_j:=t_* +  j\Bigl(  
7\vartheta + \frac{\sigma}{m}\Bigr)$ for $j=0,\dots,m$. Next, we fix
$t_* \in [t_0+\tau,t_0+2 \tau]$, $x_* \in S_D$ such that 
$$
\|v\|_{L^1(Q_-(r_0,\vartheta,t_*,x_*))}
= \max \Bigl \{\|v\|_{L^1(Q_-(r_0,\vartheta,t,x))} \::\: 
x \in S_D,\: t_0 + \tau \le t \le t_0 + 2 \tau \Bigr\}. 
$$
Since the cylinders 
$$
Q_-(r_0,\vartheta,t_0+\tau+l \vartheta,x), \qquad l
\in \N \cup \{0\},\:l \le \frac{\tau}{\vartheta},\: x \in S_D
$$
cover $[t_0 +\tau, t_0 +2 \tau] \times D$, we have 
\begin{align}
[v]_{L^1([t_0 +\tau, t_0 +2 \tau] \times D)} &= \frac{1}{\tau
  |D|}\|v\|_{L^1([t_0 +\tau, t_0 +2 \tau] \times D))}\le \frac{(n+1)(\frac{\tau}{\vartheta}+1)}{\tau
  \eps}\|v\|_{L^1(Q_-(r_0,\vartheta,t_*,x_*))} \nonumber\\
&= \frac{(n+1)(\frac{1}{\vartheta}+\frac{1}{\tau})}{
  \eps} \vartheta
|B_{r_0}(0)|\,[v]_{L^1(Q_-(r_0,\vartheta,t_*,x_*))} 
\nonumber \\
&\le \kappa_1 [v]_{L^1(Q_-(r_0,\vartheta,t_*,x_*))} \qquad
\text{with}\quad \kappa_1:= \frac{2 (n+1) |B_{r_0}(0)|}{\eps} \label{eq:36}
\end{align}
We now consider $t^* \in [t_0+3 \tau,t_0+4 \tau]$, $x \in D$
arbitrary. Then we choose $x^* \in S_D$ such that $x \in
B_{r_0}(x^*)$, and we choose $s_j$, $j=0,\dots,m$ with the
properties (\ref{eq:27}) and (\ref{eq:28}). Moreover, we fix a
sequence of points $x_j \in S_D$, $j=0,\dots,m$ such that
(\ref{eq:29}) holds with $m$ in place of $n$. This may be done, since
$m \ge n$, by repeating some of the points in the chain if necessary. 
We now define 
$$
Q_j:=Q(r_0,\vartheta,s_j,x_j) \quad \text{and}\quad 
Q_j^\pm:=Q^\pm(r_0,\vartheta,s_j,x_j)\qquad \text{for $j=0,\dots,m$.}
$$
We note that, by (\ref{eq:29}) and (\ref{eq:27}), we have
$$
|Q_j^+ \cap Q_{j+1}^-| \ge \frac{\mu \vartheta}{2} \qquad \text{for
  $j=0,\dots,m-1$.}
$$
Hence we may estimate, using Theorem~\ref{t:harnack} and the fact
that $Q_j \subset T \times U$ for $j=0,\dots,m$,  
\begin{align*}
c_1 [v]_{L^1(Q^-_j)}  &\le
\inf_{Q^+_j}v +c_2
\|g\|_{L^\infty(Q_j)}\le [v]_{L^1(Q^+_j \cap Q_{j+1}^-)} + c_2
\|g\|_{L^\infty(T \times U)}\\
&\le  \frac{|Q_{j+1}^-|}{|Q^+_j \cap Q_{j+1}^-| }[v]_{L^1(Q^-_{j+1})} + c_2
\|g\|_{L^\infty(T \times U)}\\
&\le \frac{2 |B_{r_0}(0)|}{\mu} [v]_{L^1(Q^-_{j+1})} +c_2
\|g\|_{L^\infty(T \times U)}.
\end{align*}
Iterating this estimate $m$ times and using Theorem~\ref{t:harnack}
once more, we obtain
\begin{align*}
[v]_{L^1(Q^-_0)} &\le \Bigl( \frac{2 |B_{r_0}(0)|}{c_1
  \mu}\Bigr)^{m}[v]_{L^1(Q^-_{m})}+\frac{c_2}{c_1} \sum_{k=0}^{m-1}\Bigl( \frac{2 |B_{r_0}(0)|}{c_1
  \mu}\Bigr)^{k} \|g\|_{L^\infty(T \times U)} \\
&\le \Bigl( \frac{2 |B_{r_0}(0)|}{
  \mu}\Bigr)^{m} c_1^{-(m+1)}  \inf_{Q^+_{m}}v+ \frac{c_2}{c_1} \sum_{k=0}^{m}\Bigl( \frac{2 |B_{r_0}(0)|}{c_1
  \mu}\Bigr)^{k}\|g\|_{L^\infty(T \times U)}.
\end{align*}
Hence, since $(t^*,x^*) \in Q^+_{m}$, we conclude by (\ref{eq:36}) that  
$$
v(t^*,x^*) \ge \inf_{Q^+_m} v \ge \hat c_1
[v]_{L^1(Q^-_0}-\tilde c_2\|g\|_{L^\infty(T \times U)} 
\ge \frac{\hat c_1}{\kappa_1}[v]_{L^1([t_0 +\tau, t_0 +2
  \tau] \times D)}-\tilde c_2\|g\|_{L^\infty(T \times U)}
$$
with 
$$
\hat c_1 = \Bigl( \frac{2 |B_{r_0}(0)|}{
  \mu}\Bigr)^{-m} c_1^{m+1}  \qquad \text{and}\qquad \tilde c_2 = \hat
c_1\: \frac{c_2}{c_1} \sum_{k=0}^{m}\Bigl( \frac{2 |B_{r_0}(0)|}{c_1
  \mu}\Bigr)^{k}
$$
Hence the claim follows with $\tilde c_1 = \frac{\hat c_1}{\kappa_1}$ and $\tilde c_2$ as above. Note that $\tilde c_1$
and $\tilde c_2$ only depend -- via $n$, $m$, $\mu$, $c_1$, $c_2$ and
$\kappa_1$ -- on the given quantities
$N,s,r_0,R,\eps,\tau, C_1$ and $C_2$.      
\end{proof}

The main goal of this subsection is to deduce the following Harnack
inequality for entire antisymmetric supersolutions of (\ref{eq:3-0}).

\begin{theorem}\label{harnack-anti}
 Let $r_0 \in (0,1]$, $c_\infty, R,\tau,\eps>0$ be given. Then there exist
 positive constants $K_i>0$, $i=1,2$ depending on $N,s,r_0,c_\infty,\epsilon,R,\tau$ with the following property:\\
If $D\subset\subset U\subset H$ is a pair of bounded domains with
$\dist(D,\partial
 U)\geq 4r_0$, $\diam(D) \le R$, $|D| \ge \eps$, and $v$ is entire antisymmetric supersolution of (\ref{eq:3-0}) on
 $T \times U$ with $T=[t_0,t_0 + 4\tau]$ for some $t_0 \in \R$
 and  $c \in L^\infty(T \times U)$ with 
 $\|c\|_{L^\infty} \le c_\infty$ such that $v(t) \in H^s(\R^N)$ for
 all $t \in T$, then 
 \begin{equation}
\inf_{(t,x)\in T_+\times D}v(t,x)\geq K_1 [v^{+}]_{L^{1}(T_- \times D
  )} -  K_2 \|v^{-}\|_{L^{\infty}(T \times  U)},
\end{equation}
where $T_+=[t_0+3 \tau,t_0 + 4 \tau]$ and
$T_-=[t_0+\tau,t_0+2 \tau]$.
\end{theorem}
 
The first step in the derivation of this result is the following lemma.

\begin{lemma}\label{reduction}
 Let $\beta>0$ be given, and put $H_\beta:= \{x \in H\::\:
 \dist(x,\partial H) > \beta\}$. Then there exists a continuous kernel function
 $k:\R^{N}\times\R^{N}\setminus \vartriangle \to [0,\infty)$ -- depending on
 $\beta$ -- with the
 following properties:
\begin{itemize}
\item[(i)] $k(x,y)=k(y,x)$ for all $x,y\in \R^{N}$, $x\neq y$;
\item[(ii)] $0\leq k(x,y)\leq c_{N,s}|x-y|^{-N-2s}$, for all $x,y\in \R^{N}$, $x\neq y$;
\item[(iii)] $k(x,y)\geq (1-5^{-N/2-s})c_{N,s}|x-y|^{-N-2s}$ for 
  $x,y\in \R^{N}$ with $0<|x-y|\leq \frac{\beta}{2}$;
\item[(iv)] For any antisymmetric $v \in H^s(\R^N)$ and any $\varphi \in
  \cH_0^s(H_\beta)$ we have
  \begin{equation}
\mathcal{E}(v,\varphi) = 
\cE_k(\tilde v, \phi)+2 \int_{H_\beta}\kappa_{H}(x)\tilde{v}(x)\varphi(x)\ dx\label{ziel}
\end{equation}
with $\kappa_H(x)$ as given in Lemma~\ref{sec:antisymm-supers-corr-1}
and $\tilde v = v 1_H \in \cH^s_0(H)$.
\end{itemize}
\end{lemma}

\begin{proof}
We may assume without loss that $H= \{x \in \R^N\::\: x_1>0\}$. 
For simplicity, we write $\bar x= Q(x)=(-x_1,x_2,\dots,x_N)$ for $x \in \R^N$. 
We consider $J(x,y)$ as defined in
Lemma~\ref{sec:antisymm-supers-corr-1}. Obviously we have
\begin{equation}
  \label{eq:9}
 0\leq J(x,y)\leq c_{N,s}|x-y|^{-N-2s}  \qquad \text{for $x,y \in
   \R^N,\:x \not=y.$}
\end{equation}
whereas, by Lemma~\ref{sec:antisymm-supers-corr-1}
\begin{equation}
  \label{eq:8}
J(x,y)\geq \frac{(1-5^{-N/2-s})c_{N,s}}{|x-y|^{N+2s}}\quad \text{for $x,y\in H$ with $0<|x-y|\leq \frac{\beta}{2}$
and  $\min \{x_1,y_1\} \ge \frac{\beta}{2}$.}
\end{equation}
To define $k$ with the asserted properties, we set 
$$
g:\R^{N}\times\R^{N}\setminus \vartriangle \to\R,\qquad
g(x,y):=\left\{
\begin{aligned} 
&J(x,y)&&\qquad (x,y)\in H\times H\setminus \vartriangle, \\
&0,&&\qquad \text{ otherwise},
\end{aligned} 
\right.
$$
and 
$$
s: \R^{N}\times\R^{N} \to \R,\; s(x,y):=\left\{
\begin{aligned}
&\min\{\beta-x_1,\beta-y_1\},&& \text{if $\min\{\beta-x_{1},\beta-y_{1}\} \geq 0$}, \\
&0,&&\quad  \text{ otherwise}.
\end{aligned} \right.
$$
Finally, we set $k(x,y):=g(x+s(x,y)e_1,y+s(x,y)e_1)$ for $x,y \in
\R^N,\,x \not=y$. Then $k$ is continuous, and properties (i) and (ii)
follow directly by construction and (\ref{eq:9}). To see (iii), we
note that if $0<|x-y|\leq \frac{\beta}{2}$ then also
$|\tilde{x}-\tilde{y}|\leq \frac{\beta}{2}$, where
$\tilde{x}=x+s(x,y)e_1$ and $\tilde{y}=y+s(x,y)e_1$. Furthermore we have
that $\max\{x_1,y_1\}\geq \beta$ and
therefore $\min\{x_1,y_1\} \geq
\frac{\beta}{2}$. Consequently, 
$$
k(x,y)=g(\tilde x,\tilde y) \ge (1-5^{-N/2-s})c_{N,s}|\tilde x-\tilde y|^{-N/2-s} =(1-5^{-N-2s})c_{N,s}|x-y|^{-N-2s}  
$$
by (\ref{eq:8}). It remains to show (iv): So let $v \in H^s(\R^N)$ be
antisymmetric, and let $\varphi \in \cH_0^s(H_\beta)$. Then
Lemma~\ref{sec:antisymm-supers-corr-1} gives 
\begin{equation}
\mathcal{E}(v,\varphi)=\frac{1}{2}\int_{H}\int_{H}(\tilde v(x)-\tilde v(y))({\varphi}(x)-{\varphi}(y))J(x,y)\
dx\ dy+2\int_{H_{\beta}}\kappa_{H}(x) \tilde v(x)\varphi(x)\ dx, \label{eq:11}
\end{equation}\label{eq:101-1}
whereas, since $\varphi \equiv 0$ on $\R^N \setminus H_\beta$,  
\begin{equation}
\label{eq:101}
\int_{H}\int_{H}(\tilde v(x)-\tilde v(y))({\varphi}(x)-{\varphi}(y))J(x,y)\
dx\ dy= 
\int_{H}\int_{H_{\beta}}\dots  dx\ dy +
\int_{H_{\beta}}\int_{H\setminus H_{\beta}}\dots dx \ dy .
\end{equation}
If $x\in H_{\beta}$, then for $y \in H$ we have $s(x,y)=0$ and thus
$J(x,y)=g(x,y)=k(x,y)$, while for $y\in \R^{N}\setminus
H$ we have that $k(x,y)=0$. Hence we can rewrite the first integral of
the RHS of (\ref{eq:101}) as
\begin{align*}
 \int_{H}\int_{H_{\beta}}(\tilde v(x)-\tilde v(y))({\varphi}(x)-{\varphi}(y))&J(x,y)\ dx\ dy\\
& = \int_{\R^{N}}\int_{H_{\beta}}(\tilde v(x)-\tilde v(y))({\varphi}(x)-{\varphi}(y))k(x,y)\ dx\ dy
\end{align*}
Similarly, if $y\in H_{\beta}$, then for  $x\in H\setminus H_{\beta}$
we have $s(x,y)=0$ and thus $J(x,y)=g(x,y)=k(x,y)$, while
for $x\in \R^{N}\setminus H$ we have $k(x,y)=0$. Hence we may rewrite the second integral of
the RHS of (\ref{eq:101}) as
\begin{align*}
\int_{H_{\beta}}\int_{H\setminus H_{\beta}}(\tilde v(x)-\tilde v(y))&({\varphi}(x)-{\varphi}(y))J(x,y)\ dx dy\\
&=\int_{H_\beta}\int_{\R^{N}\setminus H_{\beta}}(\tilde v(x)-\tilde v(y))({\varphi}(x)-{\varphi}(y))k(x,y)\ dx dy \\
&=\int_{\R^{N}}\int_{\R^{N}\setminus H_{\beta}}(\tilde v(x)-\tilde v(y))({\varphi}(x)-{\varphi}(y))k(x,y)\ dx dy, 
\end{align*}
where the last equality follows again since ${\varphi}=0$ on
$\R^{N}\setminus H_{\beta}$. Combining these identities, we get
\begin{align*}
 \int_{H}\int_{H}(\tilde v(x)-\tilde
 v(y))&(\varphi(x)-\varphi(y))J(x,y)\ dx dy\\
&= \int_{\R^{N}}\int_{\R^{N}}(\tilde v(x)-\tilde v(y))(\varphi(x)-\varphi(y))k(x,y)\ dx dy, 
\end{align*}
and together with (\ref{eq:11}) it follows that 
$$
 \mathcal{E}(v,\tilde{\varphi}) = \frac{1}{2}\int_{\R^{N}}\int_{\R^{N}}(\tilde{v}(x)-\tilde{v}(y))(\varphi(x)-\varphi(y))k(x,y)\ dy +2 \int_{H_\beta}\kappa_{H}(x)\tilde{v}(x)\varphi(x)\ dx,
$$
as claimed in (\ref{ziel}).
\end{proof}
 
We may now complete the 

\begin{proof}[Proof of Theorem~\ref{harnack-anti}]
Put $\beta=2 r_0$, $U_0=\{x \in
U\::\: \dist(x,D)< \beta\}\subset \subset U$, and let $k$ be the
function given by Lemma~\ref{reduction} for this choice of $\beta$. Let $v$ be an antisymmetric supersolution of (\ref{eq:3-0}) on
 $T \times U$, and consider 
$$
\tilde v: T \times \R^N \to \R,\qquad \tilde v(t,x)= \left \{
  \begin{aligned}
 &v(t,x),&&\quad (t,x) \in T \times H\\
 &0,&&\quad (t,x) \not\in T \times H.
  \end{aligned}
\right.
$$
Since $U_0 \subset H_\beta$, Lemma~\ref{reduction}(iv)
 implies that 
$$
\cE_k(\tilde v(t),\phi) \ge \int_{U_0} \Bigl([c(t,x)-2 \kappa_{H}(x)]\tilde
v(t)-\partial_t \tilde v(t)\Bigr)\phi\,dx \quad \text{for $\phi
 \in \cH_0^s(U_0)$, $\phi \ge 0$, $t \in T$,}
$$
where $0 \le \kappa_{H}(x) \le  \frac{ 4^s
  \Gamma(\frac{1}{2}+s)}{\sqrt{\pi} \Gamma(1-s)}\beta^{-2s}$ for $x
\in H_\beta$ by Lemma~\ref{sec:antisymm-supers-corr-1}. Let $d:=2 \frac{ 4^s \Gamma(\frac{1}{2}+s)}{\sqrt{\pi} \Gamma(1-s)}\beta^{-2s}+c_\infty$ and $\sigma:=\|v^{-}\|_{L^{\infty}(T \times U)}$,
and define
$w(t,x):=e^{d(t-t_0)}[\tilde v(t,x)+\sigma]$ for $t \in T,\,x \in \R^N$. Setting
$w(t)=w(t,\cdot)$ as usual, we observe that $w(t) \ge 0$ on $\R^N$ for all $t
\in T$. Moreover, for any $t \in T$ and any
 nonnegative $\phi \in \cH_0^s(U_0)$ we have
 \begin{align*}
&\cE_k(w(t),\phi) =e^{d(t-t_0)}\cE_k(\tilde v(t),\phi)\\
&\ge \int_{U_0}
  \Bigl([d+ c(t,x)-2 \kappa_{H}(x)]w(t,x) -\partial_t w(t,x)
  -e^{d(t-t_0)}\sigma [c(t,x)-2 \kappa_{H}(x)] \Bigr)\phi(x)\,dx\\
& \ge \int_{U_0}\Bigl(e^{d(t-t_0)}\sigma [2 \kappa_{H}(x)-c(t,x)]-\partial_t w(t,x)\Bigr)\phi(x)\,dx. 
\end{align*}
Hence $w$ is a nonnegative supersolution of (\ref{eq:10}) on $T \times
U_0$ with
$$
g(t,x)= e^{d(t-t_0)}\sigma [2 \kappa_{H}(x)-c(t,x)].
$$
Applying Corollary~\ref{harnack-gen} with $U_0$ in place of $U$ (noting that
$\dist(D,\partial U_0)=\beta=2 r_0$) and using the properties of $k$ given by
Lemma~\ref{reduction}, we find \mbox{$c_i=c_i(N,s,r_0,R,\eps,\tau)>0$}
such that
$$
 \inf_{T_+ \times D}w(t,x) \geq c_1 [w]_{L^{1}(T_- \times D)}
   -c_2 \|g\|_{L^\infty(T \times U_0)}
$$
We note furthermore that $[w]_{L^{1}(T_- \times D)}
\ge [v+\sigma]_{L^{1}(T_- \times D)} \ge [v^+]_{L^{1}(T_- \times D)}$
and 
$$
\inf_{T_+ \times D}w \:\le\: e^{4 \tau d}\Bigl(\inf_{T_+ \times
D}v + \sigma \Bigr),  
$$
so that 
$$
\inf_{T_+ \times D}v \ge  c_1 e^{-4 \tau d} [v^+]_{L^{1}(T_+ \times
D)}-e^{-4 \tau
  d} c_2 \|g\|_{L^\infty(T \times U_0)} -\sigma
$$
Noting furthermore that $\|g\|_{L^\infty(T \times U_0)} \le e^{4\tau d} \sigma  d,$
 we conclude that
$$
\inf_{T_+ \times D}v \ge c_1 e^{-4 \tau d} [v^+]_{L^{1}(T_+ \times
D)}-(c_2 d+1) \sigma. 
$$
Hence the assertion follows with $K_1= c_1 e^{-4 \tau d}$ and
$K_2=c_2 d+1$. Note that both constants only depend on
$N,s,r_0$, $c_\infty,\epsilon,R$ and $\tau$. 
\qedhere
\end{proof}

\subsection{A lower bound based on a subsolution estimate}\label{sub}
\label{sec:lower-boundes-via}

 The aim of this subsection is to prove the following result. 
\begin{proposition}
\label{sec:lower-bounds-via}
Let $\rho>0$, and let $\Psi$ denote the unique  positive
eigenfunction of the problem
\begin{equation}
  \label{eq:17}
\qquad \left\{ \begin{array}{rcll}
             -\Delta \Psi &=&\lambda_{1}\Psi & \text{ in } B_\rho(0), \\
                                             \Psi&=& 0 &\text{on } \partial B_\rho(0),\\
            \end{array}
\right.
\end{equation}
corresponding to the first eigenvalue $\lambda_1>0$ with $\|\Psi\|_{L^\infty(B_\rho(0))}=1$. Moreover, let
$c_\infty>0$. Then there exist
$\gamma=\gamma(N,s,\rho,c_\infty)>0$ and $q=q(N,s,\rho,c_\infty)>0$ with the
following property. 
If $T:=[t_0,t_1] \subset \R$, $x_0 \in H$ with
$\dist(x_{0},\partial H) \ge 2 \rho$, $\sigma_0>0$, $\sigma_1 \ge q \sigma_0$ and an antisymmetric
supersolution $v$ of (\ref{eq:3-0}) on $T \times B_{\rho}(x_0)$ with
$\|c\|_{L^\infty(T \times B_{\rho}(x_0))} \le c_\infty$ are
given such that 
\begin{itemize}
\item[(i)] $v(t) \in H^s(\R^N)$ for $t \in T$,
\item[(ii)] $v$ is nonnegative in $T \times B_{2\rho}(x_0)$,
\item[(iii)] $\|v^-(t)\|_{L^\infty(H \setminus B_{2\rho}(x_0))} \le \sigma_0
  e^{-(\gamma+1)(t-t_0)}$ for $t \in T$,
\item[(iv)] $v(t_0,x) \ge \sigma_1 \Psi(x-x_0)$ for $x \in
  B_{\rho(x_0)}$,
\end{itemize}
then 
\begin{equation}
  \label{eq:21}
v(t,x) \ge \sigma_1 e^{-\gamma(t-t_0)}\Psi(x-x_0) \qquad \text{for $(t,x)
  \in T \times B_{\rho(x_0)}$.}      
\end{equation}
\end{proposition}

To show this estimate, we consider the Caffarelli-Silvestre extension
of a function $v \in H^s(\R^N)$ which was introduced in \cite{CS07}. For this we consider the usual half space
$\R^{N+1}_+:= \{(x,y) \in \R^N \times \R\::\: y >0\}$. For a domain
${U_+} \subset \R^{N+1}_+$, the weighted 
Sobolev space $H^{1}(U_+;y^{1-2s})$ is given as the set
of all functions $w \in H^1_{loc}(U_+)$ such that
$$
\int_{U_+} y^{1-2s} \bigl(|w|^2+|\nabla w|^2\bigr)\,d(x,y)<\infty.
$$
In the following, we only consider the case $U_+= U \times (0,\infty)$ for some domain $U
\subset \R^N$. Then we have a well defined continuous trace map
$\tr: H^{1}(U_+;y^{1-2s}) \to H^s(U)$, see e.g. \cite{CS}. We also recall the following integration by parts
formula. If $h \in
H^1(U_+;y^{1-2s}) \cap C({\overline U} \times (0,\infty) )$ and $\tilde w \in
H^1(U_+;y^{1-2s})\cap C^2(U_+)\cap C^1({\overline U} \times (0,\infty))$ are
such that $h \equiv 0$ on $\partial U \times (0,\infty))$ and the limit $m(x):=
\lim \limits_{y \to 0}y^{1-2s} \partial_y \tilde w(x,y)$ exists in
uniform sense for $x
\in U$, then 
\begin{equation}
  \label{eq:31}
\int_{U_+}y^{1-2s} \nabla \tilde w \nabla h \,d(x,y) = - \int_{U}
m\, \tr(h)\,dx - \int_{U_+} [\div (y^{1-2s} \nabla \tilde w]
  h \,d(x,y).   
\end{equation}

Formally introducing the operator $L_s:= \div (y^{1-2s} \nabla)$, we say that a function $w \in
H^{1}(U_+;y^{1-2s})$ is (weakly)
{\em $L_s$-harmonic} on $\R^{N+1}_+$ if
$$
\int_{\R^{N+1}_+} y^{1-2s} \nabla w \nabla \phi\,dz=0 \qquad
\text{for all $\phi \in H^{1}(\R^{N+1}_+;y^{1-2s})$ with $\tr(\phi)=0$.}
$$
Standard elliptic regularity then shows that $w \in
C^\infty(\R^{N+1}_+)$ and that $\div (y^{1-2s} \nabla w)=0$ in $\R^{N+1}_+$
in pointwise sense.  
We finally recall that every function $v \in H^s(\R^N)$ has a $L_s$-harmonic
extension $w \in H^{1}(\R^{N+1}_+;y^{1-2s})$ given by 
\begin{equation}
  \label{eq:32}
w(x,y)=\int_{\R^N}v(z)G(x-z,y)\,dz \qquad \text{for $x \in \R^N, y>0,$}
\end{equation}
where
$G(x,y):=p_{N,s}\,y^{2s}\left(|x|^{2}+y^2\right)^{-\frac{N+2s}{2}}$ for $x \in \R^N, y>0$, where $p_{N,s}$ is a normalization constant, see e.g. \cite{CS07}. We need the following lemmas.

\begin{lemma}\label{kleiner}
Let $\rho>0$. Then there exists constants
$\tilde c_1=\tilde{c}_1(N,s,\rho)$ and
$\tilde{c}_{2}=\tilde{c}_{2}(N,s,\rho)$ such that the following
holds:\\
If $x_0 \in H$ satisfies $\dist(x_0,H) \ge 2 \rho$ and $v \in H^s(\R^N)$ is a continuous antisymmetric function such that $v \ge 0$ on $B_{2\rho}(x_0)$, then 
\begin{equation}\label{absch}
\frac{w(x,y)}{y^{2s}} \geq \tilde{c}_1
  \Bigl[\int_{\text{\tiny $B_{\rho}(x_0)$}}(v(z))^{\frac{1}{2}}\ dz\Bigr]^{2}-
  \tilde{c}_{2}\|v^{-}\|_{L^{\infty}(H \setminus B_{2\rho}(x_0))}\quad
\text{for $(x,y)\in B_{\rho}(x_0) \times (0,1]$,}
\end{equation}
where $w$ is the $L_s$-harmonic extension of $v$.
\end{lemma}

\begin{proof}
Since $v$ is antisymmetric, we have, by a simple change of variable,
\begin{equation}
  \label{eq:4}
w(x,y)= \int_{H} [G(x-z,y)-G(x-Q(z),y)]v(z)\,dz \qquad \text{for $x
  \in H$}. 
\end{equation}
For $x,z \in H$ and $y>0$ we have
\begin{equation}
  \label{eq:5}
 G(x-z,y) \ge G(x-z,y)-G(x-Q(z),y)
 =G(x-z,y)\left(1-\left(\frac{|x-z|^{2}+y^2}{|x-Q(z)|^2
       +y^2}\right)^{\frac{N+2s}{2}}\right)
\end{equation}
Moreover, for $x,z \in B_{\rho}(x_0)$ we have $|x-z|^2 \le 4 \rho^2$
and 
$|x-Q(z)|^2\geq |x-z|^2+4\rho^2$, so that 
\begin{equation}
  \label{eq:6}
G(x-z,y)-G(x-Q(z),y) \ge c_1\: G(x-z,y)\qquad \text{for $y \in (0,1]$}
\end{equation}
with 
$$
c_1=1-\left(1+\frac{4\rho^{2}}{1+4
    \rho^2}\right)^{-\frac{N+2s}{2}}=1-\left(\frac{1+4\rho^{2}}{1+8
    \rho^2}\right)^{\frac{N+2s}{2}}.
$$
Combining (\ref{eq:4}),~(\ref{eq:5}) and (\ref{eq:6}) and using that
$v \ge 0$ on $B_{2\rho}(x_0)$, we obtain the estimate
\begin{align}
  \label{eq:7}
&\frac{w(x,y)}{p_{N,s}y^{2s}}\geq c_1
\int_{B_{\rho}(x_0)}v(z)(|x-z|^2+y^2)^{-\frac{N+2s}{2}} \ dz \\
&-\|v^{-}\|_{L^{\infty}(H \setminus B_{2\rho}(x_0))} \int_{H\setminus
  B_{2\rho}(x_0)}(|x-z|^2+y^2)^{-\frac{N+2s}{2}}\ dz \nonumber\\  
 &\ge \frac{c_1 \Bigl(\int_{\text{\tiny$B_\rho(x_0)$}}v^{\frac{1}{2}}(z)\
   dz\Bigr)^{2}}{\int_{\text{\tiny$B_\rho(x_0)$}}(|x-z|^2+y^2)^{\frac{N+2s}{2}}\ dz}-
 \|v^{-}\|_{L^{\infty}(H\setminus B_{2\rho}(x_0))}\,
 N\omega_{N}\int_{\rho}^{\infty}r^{-1-2s}\ dr\nonumber \\
 &\geq \frac{c_1 \Bigl(\int_{\text{\tiny$B_\rho(x_0)$}}v^{\frac{1}{2}}(z)\
   dz\Bigr)^{2}}{\int_{\text{\tiny$B_\rho(x_0)$}}(|z|^2+1)^{\frac{N+2s}{2}}\ dz} -
 c_{2} \|v^{-}\|_{L^{\infty}(H\setminus B_{2\rho(x_0)})} \quad \text{for $x \in
   B_{\rho}(x_0)$, $y \in (0,1]$}  \nonumber
\end{align}
with $c_{2}:=\frac{N\omega_{N}}{2s}(\rho)^{-2s}$. Since also 
$$
\int_{B_{\rho}(0)}(|z|^2+1)^{\frac{N+2s}{2}}\ dz= N\omega_{N}
\int_{0}^{\rho}(r^2+1)^{\frac{N+2s}{2}}r^{N-1}\ dr
\le N\omega_{N}\bigl(\rho^2+1\bigr)^{\frac{3N+2s}{2}}
$$
for $x,z \in B_\rho$ and $y \in (0,1]$, we conclude that 
$$
\frac{w(x,y)}{p_{N,s}y^{2s}} \ge 
\frac{c_1}{\omega_{N}}\bigl(\rho^2+1\bigr)^{-\frac{3N+2s}{2}}
  \Bigl(\int_{B_{1}}v^{\frac{1}{2}}(z)\ dz\Bigr)^{2} -  \tilde c_{2} \|v^{-}\|_{L^{\infty}(H\setminus B_2)}
$$
for $x \in B_1$ and $y \in (0,1]$. Hence the claim follows with
$\tilde c_1 = \frac{c_1 p_{N,s}}{\omega_{N}}\bigl(\rho^2+1\bigr)^{-\frac{3N+2s}{2}}$ and $\tilde c_2=c_2 p_{N,s}$.
\end{proof}

\begin{lemma}
\label{sec:antisymm-supers-line-2}
Let $U \subset H$ be a bounded Lipschitz domain, $T:=(t_{0},T_{0})$,
$t_{0}<T_{0}$, $c \in L^\infty(T \times U)$, and let $v$ be a supersolution of (\ref{eq:3-0}) 
on $T \times U$ such that $v(t) \in H^s(\R^N)$ for all $t \in
T$. Moreover, let $w(t) \in H^1(\R^{N+1}_+,y^{1-2s})$ be the $L_s$-harmonic extension of
$v(t)$ given by (\ref{eq:32}) for each fixed time $t \in T$. Then for 
every nonnegative $\Phi \in H^1(\R^{N+1}_+,y^{1-2s})$ with $\phi:= \tr(\Phi) \in
\cH_0^s(U)$ and every $t \in T$ we have 
\begin{equation}
  \label{eq:13}
\int_{\R^{N+1}_+} y^{1-2s} \nabla w \nabla \Phi\,d(x,y) \ge d_{s}\int_{U} (c(t,x)v-\partial_{t}
v)\phi\,dx,
\end{equation}
where $d_{s}=2^{1-2s}\Gamma(1-s)/\Gamma(s)$.
\end{lemma}

\begin{proof}
In case $\Phi$ is the $L_s$-harmonic extension of $\phi$, we have
$$
\int_{\R^{N+1}_+} y^{1-2s} \nabla w \nabla \Phi\,d(x,y)= d_{s}\cE(w,\phi)
$$  
with $d_{s}$ as stated (see e.g. \cite{CS07} or \cite[Remark 3.11]{CS}) and therefore (\ref{eq:13})
is true. On the other hand, since $w$ is $L_s$-harmonic,
$$
\int_{\R^{N+1}_+} y^{1-2s} \nabla w \nabla \Theta\,d(x,y)=0\qquad \text{for every
$\Theta \in H^1(\R^N_+,y^{1-2s})$ with $\tr(\Theta)=0$.}
$$
Hence the assertion follows.   
\end{proof}

\begin{lemma}\label{subloesung}
Let $\rho>0$, $c_\infty>0$, and let $\Psi$ be defined as in
Proposition~\ref{sec:lower-bounds-via} w.r.t. $\rho$.  Then there exists 
$\gamma= \gamma(N,s,\rho,c_\infty)>0$ such that the following holds:\\ 
If $T=[t_0,t_1] \subset \R$, $x_0 \in H$ with $\dist(x_0,\partial
H)\ge 2\rho$ and $c \in L^\infty([t_0,t_1] \times
B_\rho(x_0))$ with $\|c\|_{L^\infty} \le c_\infty$ are
given and $v$ is an antisymmetric supersolution of (\ref{eq:3-0}) on $[t_0,t_1]
\times B_\rho(x_0)$ such that $v(t) \in H^s(\R^N)$ for $t \in
[t_0,t_1]$,  
\begin{equation}
  \label{eq:12}
v(t_0,x) \ge \sigma\, \Psi(x-x_0) \qquad \text{for $x \in B_\rho(x_0)$ with
some constant $\sigma>0$,}  
\end{equation}
and 
\begin{equation}
  \label{eq:16}
  \begin{aligned}
&\text{the $L_s$-harmonic extension $w(t)$ of $v(t)$ is nonnegative}\\
&\text{on
$B_\rho(x_0) \times [0,1]$ for all $t \in [t_0,t_1]$,}      
  \end{aligned}
\end{equation}
then 
$$
v(t,x) \ge \sigma e^{-\gamma (t-t_0)} \Psi(x-x_0) \qquad \text{for $(t,x) \in
T \times  B_\rho(x_0)$.}
$$
\end{lemma}

\begin{proof}
Without loss of generality, we may assume that $x_0=0$, $t_0=0$ and $\sigma=1$, and
we put $B_\rho=B_\rho(0)$. Let
$\lambda_1>0$ be defined by (\ref{eq:17}), and let $f: [0,\infty) \to \R$ denote the solution of the initial value problem
$$
 (D)\qquad \left\{
   \begin{aligned}
&f''+\frac{1-2s}{y}f'-\lambda_1 f=0, \\
&f(0)=1,\\
&\lim\limits_{y\to\infty}f(y)=0,     
   \end{aligned}
\right.
$$
which is uniquely given by 
$$
 f(y)=\kappa_1 y^{2s}\int_{0}^{\infty}\frac{\cos(\tau)}{(\tau^{2}+ \lambda_1 y^{2})^{\frac{1+2s}{2}}}\
 d\tau, \text{for $y\ge 0$ with $\kappa_1 =\lambda_1^s d_{s}$,}
$$
with $d_{s}$ as in Lemma \ref{sec:antisymm-supers-line-2}. We note that $f$ is a scalar multiple of a rescaled MacDonald function
(or modified Bessel function of the second kind),
see e.g. \cite{Watson}. We also note that $f$ is strictly decreasing
on $[0,\infty)$. Moreover, the
limit 
\begin{equation*}
\kappa_2:=\lim \limits_{y\to
  0} \frac{y^{1-2s}f'(y)}{1-f(1)}  \le 0
\end{equation*}
exists and only depends on $s$ and $\rho$ (via $\lambda_1$). We now
put $\gamma= c_\infty- \kappa_2 +1$, and we let 
$$
\text{$\tilde w: [0,t_1] \times \R^{N+1}_+ \to \R$ be
  defined by}\;  
\tilde w(t,x,y)=\left \{
\begin{aligned}
&e^{-\gamma t}\Psi(x)\frac{f(y)-f(1)}{1-f(1)},&&\quad x \in
B_\rho;\\
&0,&&\quad x \not \in B_\rho.  
\end{aligned}
\right.
$$
Putting $\tilde w(t)= \tilde w(t,\cdot,\cdot)$ as usual, we then have 
\begin{align}
 &L_{s}\tilde w(t) =\frac{1}{1-f(1)} e^{-\gamma t}\Bigl( y^{1-2s}\Delta_{x}\tilde w
 +(1-2s)y^{-2s}\partial_{y}\tilde w+y^{1-2s}\partial_{yy}\tilde
 w\Bigr) \label{eq:26}\\
 &=
\frac{1}{1-f(1)} e^{-\gamma t}  y^{1-2s}\,\Psi(x)\Bigl(-\lambda_1 f + \frac{1-2s}{y}f'
 +f'' + \lambda_1 f(1) \Bigr) \ge 0 \quad \text{on
   $B_\rho\times(0,\infty)$} \nonumber
\end{align}
for $t \in [0,t_1]$. Moreover, we have
\begin{align}
  \label{eq:19}
&\tilde w(0,x,0)= \Psi(x) \le v(0,x)\quad \text{for $x \in
  B_\rho$,}\\
&\lim \limits_{y\to
  0}y^{1-2s} \partial_y \tilde w(t,x,0)= \kappa_2 e^{-\gamma t }\Psi(x) = \kappa_2
\tilde w(t,x,0) \quad \text{for $x \in B_\rho$, $t \in [0,t_1]$, and}   \label{eq:25}\\
&\tilde w(t) \equiv 0 \le w(t) \qquad \text{on $\partial B_\rho \times
  (0,1] \cup B_\rho \times \{1\}$} \nonumber
\end{align}
by assumption and by construction of $\tilde w$.  In the following, we
consider 
$$
h(t) \in H^1 \in (\R^{N+1}_+,y^{1-2s}), \quad h(t,x)= 
\left\{
\begin{aligned}
&(w-\tilde w)^-(t,x,y),&& (x,y) \in B_\rho \times [0,1],\\
&0, && \text{elsewhere}.   
\end{aligned}
\right.
$$
Moreover, we will write $g$ resp. $\tilde v$ for the traces of $h$ and
$\tilde w$, respectively. Then, as a consequence of (\ref{eq:31}), (\ref{eq:13}),
(\ref{eq:26}) and (\ref{eq:25}), we have  
\begin{align*}
0 &\ge -d_{s}^{-1}\int_{\R^{N+1}_+}y^{1-2s}|\nabla h|^2 d(x,y) \\
&=d_{s}^{-1}\int_{\R^{N+1}_+}y^{1-2s}\nabla w \nabla h d(x,y)-
d_{s}^{-1}\int_{\R^{N+1}_+}y^{1-2s}\nabla \tilde w \nabla h d(x,y)\\
&\ge \int_{B_{\rho}}\Bigl([c(t,x)v-\partial_t v]g + \kappa_2 \tilde v
 g\Bigr)\,dx\\
&= \int_{B_{\rho}}\Bigl([c(t,x)(v-\tilde v) -\partial_t (v-\tilde v)]g
+ [\kappa_2+c(t,x) +\gamma]\tilde v g\Bigr)\,dx\\
&\ge -c_\infty \int_{B_\rho} g^2(t,x)\,dx + \frac{1}{2}\frac{d}{dt} \int
g^2(t,x)\,dx \qquad \text{for $t \in [0,t_1]$.} 
\end{align*}
Hence $\frac{d}{dt} \int_{B_\rho} g^2(t,x)\,dx \le 2c_\infty \int_{B_\rho}
g^2(t,x)\,dx$ for $t \in [0,t_1]$. Since furthermore $\int_{B_\rho} g^2(0,x)\,dx = 0$ as a consequence of (\ref{eq:19}), we conclude
that $\int_{B_\rho} g^2(t,x)\,dx =0$ and therefore $g(t) \equiv 0$ on
$B_\rho$ for all $t \in T$. Hence $v(t,x) \ge e^{-\gamma t} \Psi(x)$ for $(t,x) \in T
\times  B_\rho(0)$, as claimed.
\qedhere
\end{proof}

We may now complete the

\begin{proof}[Proof of Proposition~\ref{sec:lower-bounds-via}]
For given $\rho,c_\infty>0$, let $\tilde c_i$, $i=1,2$ be given by
Lemma~\ref{kleiner}, and let $\gamma$ be given by
Lemma~\ref{subloesung}. Moreover, let
$$
q= \frac{2 \tilde c_2}{\tilde{c}_1}\Bigl[\int_{\hspace{0.4cm}\vspace{-0.2cm}B_{\rho}(0)}
\Psi(z)^{\frac{1}{2}}\ dz\Bigr]^{-2}
$$
Next, let $T:=[t_0,t_1] \subset R$, $\sigma_0>0$ and $\sigma_1 \ge q \sigma_0$, and
let $v$ be an antisymmetric
supersolution of (\ref{eq:3-0}) on $T \times B_{\rho}(x_0)$ satisfying
assumptions (i)- (iii). Suppose by contradiction that   
\begin{equation}
  \label{eq:21_1}
v(t,x)= \sigma_2 e^{-\gamma(t_*-t_0)}\Psi(x_*-x_0) \quad \text{for some $\sigma_2
  \in (\frac{\sigma_1}{2},\sigma_1)$, $t_* \in T$ and $x_* \in B_{\rho}(x_0)$.}      
\end{equation}
We may assume that $t_*>t_0$ is chosen minimally with this property,
so that
$$
v(t,x) \ge \sigma_2 e^{-\gamma(t-t_0)}\Psi(x-x_0)\qquad \text{for $(t,x)
  \in [t_0,t_*] \times B_{\rho(x_0)}$.}      
$$
Let $w$ denote the $L_s$-harmonic extension of $v$ given by
(\ref{eq:4}) for each fixed time $t \in T$. Then Lemma~\ref{kleiner} implies that 
\begin{align*}
y^{-2s} w(t,x,y)&\geq \tilde{c}_1
  \Bigl[\!\!\int_{\hspace{0.4cm}\vspace{-0.2cm} B_{\rho}(x_0)}(v(t,z))^{\frac{1}{2}}\ dz\Bigr]^{2}-
  \tilde{c}_{2}\|v^{-}(t)\|_{L^{\infty}(H \setminus
    B_{2\rho}(x_0))}\\ 
&\geq \tilde{c}_1 \sigma_2 e^{-\gamma(t-t_0)}
  \Bigl[\!\! \int_{\hspace{0.4cm}\vspace{-0.2cm} B_{\rho}(0)}\Psi(z)^{\frac{1}{2}}\ dz\Bigr]^{2}-
  \tilde{c}_{2}\sigma_0 e^{-(\gamma+1)(t-t_0)}\\ 
&\geq  \sigma_0 e^{-\gamma(t-t_0)} \Bigl( \frac{q \tilde{c}_1}{2} \Bigl[\!\!\int_{\hspace{0.4cm}\vspace{-0.2cm} B_{\rho}(0)}\Psi(z)^{\frac{1}{2}}\ dz\Bigr]^{2}-
  \tilde{c}_{2}\Bigr) \ge 0
\end{align*}
for $t \in [t_0,t_*]$ and $(x,y)\in B_{\rho}(x_0) \times
  (0,1]$. Hence, by Lemma~\ref{subloesung}, 
$$
v(t,x) \ge \sigma_1 e^{-\gamma (t-t_0)} \Psi(x-x_0) \qquad \text{for $(t,x)
  \in [t_0,t_*] \times  B_\rho(x_0)$.}
$$
This contradicts (\ref{eq:21_1}), and thus the proof is finished.
\end{proof}

\subsection{\large{\textit{Main Result on entire antisymmetric supersolutions}}}\label{mainlinear}

We recall from Section~\ref{sec:notation} that, for a subset $D
\subset \R^N$, $\inrad(D)$ denote the supremum of
all $r>0$ such that every connected component of $D$ contains a ball 
$B_{r}(x_{0})$ with $x_0 \in D$. Note that $\inrad(D) \ge \rho>0$
implies that every connected
component of $D$ has at least measure $|B_\rho(0)|$, so in this case $D$ has only
finitely many connected components if it has finite measure.  
 
\begin{theorem}\label{sec:haupt}
Let $\rho>0$ and $c_\infty>0$ be given. Moreover, let
$\gamma=\gamma(N,s,\rho,c_\infty)>0$, $q=q(N,s,\rho,c_\infty)>0$ be as in
Proposition~\ref{sec:lower-bounds-via},  and let $\delta>0$ be such that the conclusions
of Proposition \ref{sec:klein1-1} hold with $\gamma+1$  in place of $\gamma$.\\ 
Then for any $\tau,r_0,R>0$ there exists $\mu>0$ such that the following holds:\\ 
If $D\subset\subset U \subset H$ are bounded open sets with
$|U|<\infty$, $\inrad(D)>2\rho$, $\diam(D)\le R$,
$|U\setminus \overline{D}|< \delta$ and $\dist(D,\partial
U)>4r_0$ and if $v$ is an entire antisymmetric supersolution of
(\ref{eq:3-0}) on $[t_0,\infty) \times U$ for some $t_{0} \in \R$ with
$\|c\|_{L^{\infty}([t_0,\infty) \times U)} \leq c_\infty$ such
that $v(t) \in H^s(\R^N)$ for $t \in [t_0,\infty)$, $v$ is nonnegative on
$[t_{0},t_{0}+8\tau]\times \overline{D}$ and 
\begin{equation}\label{sec:bedingung}
\|v^{-}(t_{0},\cdot)\|_{L^{\infty}(U\setminus \overline{D})} \leq \mu
[v]_{L^{1}((t_{0}+\tau,t_{0}+2\tau)\times D_*)}, \text{ for each connected component $D_*$ of $D$},
\end{equation}
 then:
\begin{enumerate}
\item[(i)] $v(t,x)>0$ for all $(t,x)\in [t_{0},\infty)\times \overline{D}$
\item[(ii)] $\|v^{-}(t)\|_{L^{\infty}(U)}\to 0$ for $t\to\infty$.
\end{enumerate}
\end{theorem}

\begin{proof}
We let $\rho,\gamma,q,\delta,\tau,R$ be given with the
properties stated in the theorem, and we put $\eps=|B_{2\rho}(0)|$. Let $K_1,K_2$ -- depending on these
quantities -- be given as in Theorem~\ref{harnack-anti}. We fix
$\mu>0$ sufficiently small such that 
\begin{equation}
  \label{eq:24}
\Bigl(\frac{K_1 }{\mu} -K_2 \Bigr)>q\qquad \text{and}\qquad 
\frac{K_1 |B_\rho(0)|}{(2R)^N}  \Bigl(\frac{K_1}{\mu}
-K_2  \Bigr) [\Psi]_{L^1(B_\rho(0))} - K_2>0,
\end{equation}
where $\Psi$ is given in Proposition~\ref{sec:lower-bounds-via} depending on
$\rho$. Next, we consider $D\subset\subset U \subset H$ and an antisymmetric supersolution $v$ of
(\ref{eq:3-0}) on $[t_0,\infty) \times U$ with the properties stated in
the theorem, which implies in particular that $\eps \le |D_*| \le
(2R)^N$ for every connected component $D_*$ of $D$. We put $\sigma_0= \|v^{-}(t_0)\|_{L^\infty (U \setminus
  \overline D)}$ and
$$
T_0:= \sup \{t \ge t_0 + 8 \tau\::\: \text{$v > 0$ in $[t_0,t)
  \times \overline D$}\},
$$
so that $t_0 + 8 \tau \le T_0 \le \infty$ by assumption. Applying Proposition \ref{sec:klein1-1}, we get
\begin{equation}
\|v^{-}(t)\|_{L^{\infty}(U)}=\|v^{-}(t)\|_{L^{\infty}(U\setminus \overline{D})}
\leq \sigma_{0} e^{-(\gamma+1)(t-t_{0})}\qquad \text{ for all } t\in[t_{0},T_{0}).\label{sec:bringts}
\end{equation}
To prove (i), we suppose by contradiction that
$T_0<\infty$. Then there exists a connected component $D_*$ of $D$ and 
$x_* \in \overline D_*$ such that  
\begin{equation}
  \label{eq:20}
\text{$v > 0$ in $[t_0,T_0) \times \overline D_*$}\qquad \text{and}\qquad
\text{$v(T_0,x_*)=0$.}
\end{equation}
Let $U_*$ be the connected component of $U$ with $D_* \subset U_*$. 
Since $v\geq 0$ on $[t_{0},t_{0}+8\tau)\times \overline{D_*}$, we have, by Theorem
\ref{harnack-anti}, (\ref{sec:bedingung}) and Proposition~\ref{sec:klein1-1},
\begin{align}
&\inf_{[t_{0}+3\tau,t_{0}+4\tau]\times \overline{D_*}}v \geq K_1
[v^{+}]_{L^{1}([t_{0}+\tau,t_{0}+2\tau]\times D_*)} -K_2 
\|v^{-}\|_{L^{\infty}([t_0,t_{0}+4\tau] \times U_*)} \nonumber\\
&\ge K_1
[v]_{L^{1}([t_{0}+\tau,t_{0}+2\tau]\times D_*)} -K_2 
\|v^{-}\|_{L^{\infty}([t_0,t_{0}+4\tau] \times [U \setminus \overline
  D])} \nonumber\\
&\ge \frac{K_1}{\mu} \|v^{-}(t_{0},\cdot)\|_{L^{\infty}(U\setminus
  \overline{D})}
-K_2 \|v^{-}(t_{0},\cdot)\|_{L^{\infty}(U \setminus
  \overline{D})}
= \Bigl(\frac{K_1 }{\mu} -K_2 \Bigr)\sigma_0=:\sigma_1. \label{eq:18}
\end{align}
We fix $x_{0} \in D_*$ such that $B_{2 \rho}(x_0) \subset D_*$,
which is possible by assumption. Since $\sigma_1 \ge q \sigma_0$ by
(\ref{eq:24}), the estimates (\ref{eq:18}) and
(\ref{sec:bringts}) allow us to apply
Proposition~\ref{subloesung} with $t_0+4 \tau$ in place of $t_0$, which yields
\begin{equation}
  \label{eq:23}
v(t,x) \ge \sigma_1 e^{-\gamma(t-t_0-4 \tau)}\Psi(x-x_0) \qquad \text{for every
  $x \in B_{\rho}(x_0)$, $t \in [t_{0}+4\tau, T_0]$.}  
\end{equation}
With the help of Theorem \ref{harnack-anti},$\,$ (\ref{sec:bringts})
and (\ref{eq:23}), we find that 
\begin{align*}
&v(T_0,x_*)  \geq K_1
[v]_{L^{1} ([T_0 - 3 \tau,T_0-2 \tau] \times D_*)} -K_2 
\|v^{-}\|_{L^{\infty}([T_0-4\tau,T_0] \times U_*)} \nonumber\\
&\ge K_1 \sigma_1 e^{-\gamma (T_0-6\tau-t_{0})}\frac{|B_\rho(0)|}{|D_*|}[\Psi]_{L^1(B_\rho(0))} -
  K_2 \sigma_0 e^{-(\gamma+1)(T_0-4 \tau-t_0)} \\
&\ge \sigma_0 e^{-\gamma(T_0-4 \tau - t_{0})}\Bigl[ \frac{K_1 |B_\rho(0)|}{(2R)^N}  \Bigl(\frac{K_1}{\mu}
-K_2 \Bigr) [\Psi]_{L^1(B_\rho(0))} - K_2 \Bigr] 
>0,
\end{align*} 
by our choice of $\mu$ in (\ref{eq:24}), contradicting
(\ref{eq:20}). We conclude that $T_0=\infty$. In particular, (i)
holds, and (ii) follows since, by (\ref{sec:bringts}), 
\begin{equation}
\|v^{-}(t)\|_{L^{\infty}(U)}=\|v^{-}(t)\|_{L^{\infty}(U\setminus \overline{D})}
\leq \sigma_{0} e^{-(\gamma+1)(t-t_{0})}\qquad \text{ for all } t\in[t_{0},\infty).\label{sec:bringts-1}
\end{equation}
\qedhere
\end{proof}

\section{Proof of the main symmetry result}
\label{sec:proof-main-symmetry}

In this section we complete the proof of Theorem~\ref{sec:goal}. With the tools developed in Section~\ref{1}, we may
follow the main lines of the moving plane method as developed by
Pol\'a\v cik in \cite{P07}, but some steps in the
argument -- in particular the proofs of Lemma~\ref{sec:schritt2} and Proposition~\ref{sec:proof-main-symmetry-2} below -- differ significantly
from \cite{P07}. This is due to the fact that, contrary to \cite{P07},
we do not a priori assume the existence of an element $\phi \in \omega(u)$ with
$\phi>0$. For $\lambda \in \R$, we use the notations  
$$
\Omega_\lambda=
\{x \in \Omega\::\: x_1 >\lambda\},\quad  H_\lambda:=\{x \in \R^N\::\: x_1> \lambda\},\quad
T_\lambda= \partial H_\lambda \quad \text{and} \quad \Gamma_\lambda=
T_\lambda \cap \Omega.
$$
Moreover, we let $Q_\lambda: \R^N \to \R^N$ denote the reflection at
$T_\lambda$ given by $Q_\lambda(x)= (2\lambda-x_1,x_2,\dots,x_N)$. For
a function $z: \R^N \to \R$, we put 
$$
z^\lambda= z \circ Q_\lambda\::\: \R^N \to \R
$$
and 
$$
V_\lambda z: \R^N \to \R,\qquad V_\lambda z(x) = z^\lambda(x) - z(x).
$$
We now assume that the hypotheses $(D1)$ and $(F1),(F2)$ are
satisfied, and we let $u$ be a nonnegative global solution of $(P)$ satisfying $(U_1)$
and $(U_2)$. We set 
$$
l:= \max \{x_1\::\: (x_1,x') \in \Omega \quad \text{for some $x' \in
  \R^{N-1}$}\},
$$ 
and we fix $\lambda \in [0,l)$ for the moment. As discussed in
Remark~\ref{sec:antisymm-supers-corr-4}, the function $v:=V_\lambda u$ is an
entire antisymmetric supersolution of the problem  
\begin{equation}
  \label{eq:3}
\partial_t v +\hl v = c_\lambda(t,x) v 
\end{equation}
in $(0,\infty) \times \Omega_\lambda$ with   
$$
c_\lambda(x,t)= \left \{
  \begin{aligned}
  &\frac{f(t,x,u^\lambda(x))-f(t,x,u(x))}{u^\lambda(x)-u(x)},&&\qquad u^\lambda(t,x) \not= u(t,x);\\
  &0, &&\qquad u^\lambda(t,x)= u(t,x). 
  \end{aligned}
\right.
$$
Here the term entire antisymmetric supersolution refers to the notion
defined in the beginning of Section~\ref{1} with respect to the half
space $H=H_\lambda$. Indeed, for $\lambda \in [0,l)$ and this choice
of $H$, (\ref{eq:39}) and (\ref{eq:22}) are satisfied as a consequence
of assumptions $(D1)$ and $(F2)$. Moreover, as a consequence of $(F1)$
and $(U1)$, there exists $c_\infty>0$ such that 
$$
\|c_\lambda\|_{L^\infty((0,\infty) \times \Omega_\lambda)} \le
c_\infty \qquad \text{for every $\lambda \in [0,l)$.}
$$
In the following, we fix $c_\infty$ with this property. We also note
that $[V_\lambda u](t) \in H^s(\R^N)$ for all $t \in (0,\infty)$. For $\lambda\in [0,l)$, we now consider the following statement:
\begin{equation*}
 (S_{\lambda})\qquad \|(V_{\lambda}u)^{-}(t)\|_{L^{\infty}(H_{\lambda})}\to 0 \text{ as } t\to \infty.
\end{equation*}
Our aim is to show, via the method of moving planes, that $(S_\lambda)$ holds for every $\lambda \in
[0,l)$. We need the following lemmas.

\begin{lemma}\label{sec:schritt1}
 There is $\delta>0$ such that for each $\lambda \in [0,l)$ the following statement holds. If $K$ is a closed subset of $\Omega_{\lambda}$ with $|\Omega_{\lambda}\setminus K|<\delta$ and there is $t_{0}\geq 0$ such that $V_{\lambda}u(t)\geq 0$ on $K$ for all $t\geq t_{0}$, then 
\begin{equation}
  \|(V_{\lambda}u)^{-}(t)\|_{L^{\infty}(H_{\lambda})}\leq e^{-(t-t_{0})} \|(V_{\lambda}u)^{-}(t_{0})\|_{L^{\infty}(H_{\lambda})},
\end{equation}
for all $t\geq t_{0}$. In particular $(S_{\lambda})$ holds if $\lambda<l$ is sufficiently close to $l$.
\end{lemma}

\begin{proof}
This follows immediately by applying Proposition~\ref{sec:klein1-1} to
$\gamma=1$, $c_\infty>0$ as fixed above, $H=H_\lambda$ and $U=\Omega_{\lambda}\setminus K$. Note that
the number $\delta>0$ given by Proposition~\ref{sec:klein1-1} in this
case does not depend on $\lambda$ and $K$. The second statement of the
lemma follows since $|\Omega_{\lambda}|<\delta$ if $\lambda$ is close to $l$.
\qedhere
\end{proof}

\begin{lemma}\label{sec:schritt2}
 Suppose $\lambda_{0}\in[0,l)$ is such that $(S_{\lambda})$ holds for
 all $\lambda\in (\lambda_{0},l)$. Then we have: 
 \begin{itemize}
 \item[(i)] $(S_{\lambda_{0}})$ holds.
 \item[(ii)] For each $z\in \omega(u)$ we have either $V_{\lambda_{0}}z>0$ on $\Omega_\lambda$ or
$V_{\lambda_{0}}z \equiv 0$ on $\R^N$. 
\item[(iii)] If $\lambda_0>0$, then for each $z\in \omega(u)$ we have
  either $V_{\lambda_{0}}z>0$ on 
$\Omega_\lambda$ or $z \equiv 0$ on $\R^N$. 
\end{itemize}
\end{lemma}

\begin{proof}
(i) Since the set $\{u(t)\::\:t \ge 0\}$ is relatively compact in
 $C(\overline{\Omega})$, the statement $(S_{\lambda})$ is equivalent to 
$V_{\lambda}z\geq 0 \text{ on } \Gamma_{\lambda}  \text{ for all }
z\in \omega(u).$ Hence $(S_{\lambda_{0}})$ holds by assumption and continuity of all $z
\in \omega(u)$.\\
(ii) {\em Step one:} We first claim that on each connected component $U$ of
 $\Omega_{\lambda_0}$ we either have $V_{\lambda_{0}}z>0$ on $U$ or
 $V_{\lambda_{0}}z\equiv 0$ on $U$. To prove this, we fix $z\in \omega(u)$ and a connected component $U$ of
$\Omega_{\lambda_0}$ such that $V_{\lambda_{0}}z\not\equiv 0$ on
$\Omega_{\lambda_{0}}$. Since $V_{\lambda_{0}}z\geq 0$, there exists
$x_0 \in U$ and $\rho>0$ such that $B:= B_\rho(x_0)
\subset\subset \Omega_{\lambda_{0}}$ and $V_{\lambda_{0}}z>0$ on
$\overline{B}$. Since $z\in \omega(u)$, there exists a sequence of numbers $t_n>0$,
such that $t_{n}\to \infty$ and $u(t_{n})\to z$ in
$C(\overline{\Omega})$, hence also $V_{\lambda_{0}}u(t_{n})\to
 V_{\lambda_{0}}z$ in $C(\overline{\Omega_{\lambda_0}})$ as $n \to \infty$. Consequently, there exists $\sigma>0$ and
 $n_{0} \in \N$ such that 
\[
 V_{\lambda_{0}}u(t_{n},x)>2\sigma \qquad \text{for $x\in \overline{B},\, n>n_{0}.$}
\]
By the equicontinuity property $(U2)$, there exists $\tau>0$ such that 
\begin{equation}
  \label{eq:40}
 V_{\lambda_{0}}u(t,x)>\sigma
\qquad \text{for $x\in \overline{B},\,t\in[t_{n}-4\tau,t_{n}] ,\,n>n_{0}.$}
\end{equation}
Now fix a subdomain $D \subset \subset U$. Applying Proposition~\ref{harnack-anti} with $U=\Omega_{\lambda_0}$,
$t_0=t_{n}-4 \tau$ and using (\ref{eq:40}), we get 
\begin{align*}
\inf_{x \in D}V_{\lambda_0}u(t_n,x)&\geq K_1 [(V_{\lambda_0}u)^{+}]_{L^{1}([t_{n}-4\tau,t_{n}-3 \tau]\times D
  )} -  K_2 \sup_{t\in T}\|(V_{\lambda_0}u)^{-}(t,\cdot)\|_{L^{\infty}(U)}\\
&\ge K_1 \sigma \frac{|B|}{|D|} -K_2 \|v^{-}(t,\cdot)\|_{L^{\infty}(T
  \times U)} \qquad \text{for $n>n_0$}
\end{align*}
with suitable constants $K_1,K_2>0$ independent of $n$. Since $(S_{\lambda_{0}})$ holds, we conclude that  
$$
\inf_{x \in D}V_{\lambda_0}z= \lim_{n \to \infty}\: \inf_{x \in D}V_{\lambda_0}u(t_n,x) \ge K_1
\sigma \frac{|B|}{|D|}>0.
$$
Since $D \subset \subset U$ was chosen arbitrarily, we conclude that
$V_{\lambda_0}z>0$ in $U$. This shows the claim.\\
{\em Step two:} Let 
$z\in \omega(u)$ be such that
$$
U_z := \{x \in \Omega_{\lambda_{0}}\::\: [V_{\lambda_{0}}z](x)=0\}
$$
is nonempty. To finish the proof of (ii), we need to show that $V_{\lambda_{0}}z\equiv 0$ on
$\R^N$. We suppose by contradiction that this is false; then
there exists a compact set $\cK \subset H_{\lambda_0} \setminus \overline{U_z}$ of positive
measure such that 
\begin{equation}
  \label{eq:63}
\inf_{\cK}V_{\lambda_0}z>0.   
\end{equation}
By Step one above, $U_z$ is an open set. Hence we may fix a nonnegative function $\varphi\in
C^{\infty}_{c}(U_z)$, $\varphi \not \equiv 0$, and we set $D:= \supp\, \phi$. Moreover, we fix $\rho>0$ with $\dist(D, \partial U_z)> 2\rho$, and we note that there exists $M>0$ such that 
\begin{equation}
  \label{eq:42}
\Bigl| \int_{B_{\rho}(x)}\frac{\varphi(x)-\varphi(y)}{|x-y|^{N+2s}}\
dy \Bigr| <M
 \qquad \text{for all $x \in \R^N$,}  
\end{equation}
see e.g. \cite[Lemma 3.5]{NPV11}). In the following, we put $v=
V_{\lambda_0} u$ and $H=H_{\lambda_0}$. Moreover, we consider $J$ and $\kappa$ as defined in
Lemma~\ref{sec:antisymm-supers-corr-1} for this choice of $H$.  By Lemma \ref{sec:antisymm-supers-corr-1} we have
\begin{align}
  \label{eq:61}
 \cE(v(t),\varphi)&=\frac{1}{2}\int_{H}\int_{H}(v(t,x)-v(t,y))(\varphi(x)-\varphi(y))J(x,y)\
 dxdy\\
&+ 2\int_{H}v(t,x)\kappa_H(x)\varphi(x)\ dx, \nonumber
\end{align}
where 
$$
\int_{H}v(t,x)\kappa_H(x)\varphi(x)\ dx \le
\kappa_s\|\phi\|_{L^1(U_z)} \|v(t)\|_{L^\infty(U_z)} \quad
\text{with}\;\kappa_s:=\frac{ 4^s
  \Gamma(\frac{1}{2}+s) }{\sqrt{\pi} \Gamma(1-s)}
(2\rho)^{-2s}.
$$
To estimate the double integral on the right hand side of (\ref{eq:61}), we
put 
$$
\cH_1:= \{(x,y) \in H \times H\::\: |x-y|\le \delta\}, \qquad \cH_2:= H
\times H \setminus \cH_1
$$ 
and $D_\rho:= \{x \in \R^N\::\: \dist(x,D)\le \rho\}$. Then  
\begin{align*}
 \int_{\cH_1} & (v(t,x)-v(t,y))(\varphi(x)-\varphi(y))J(x,y) dxdy\\
&=  \!\!\int_{{\substack{|x-y|\leq\delta, \\ x,y\in D_{\rho}}}}\!\! (v(t,x)-v(t,y))(\varphi(x)-\varphi(y))J(x,y) dxdy\\
&=c_{N,s}\biggl(2\int_{D_{\rho}}v(t,x)\int_{B_{\rho}(x)}\frac{\varphi(x)-\varphi(y)}{|x-y|^{N+2s}}\
  dydx \\
&\qquad\qquad\qquad\qquad-\int_{\substack{|x-y|\leq\delta, \\ x,y\in D_{\rho}}}\frac{(v(t,x)-v(t,y))(\varphi(x)-\varphi(y))}{|x-Q_{\lambda_{0}}(y)|^{N+2s}}\ dxdy\biggr)\\
&\le 2 c_{N,s} M  \: |D_\rho| \:\|v(t)\|_{L^\infty(U_z)} +
\frac{4|D_\rho|^2}{(2\rho)^{N+2s}} \, \|\phi\|_{L^\infty(U_z)} \,\|v(t)\|_{L^\infty(U_z)},
\end{align*}
where we used the fact that $|x-Q_{\lambda_{0}}(y)| \ge 2 \rho$ for
every $x,y \in D_\rho$. To estimate the integral over $\cH_2$, we
first note that 
$$
\sup_{x \in H}\; \int_{H \setminus
  B_{\rho}(x)}J(x,y)\ dydx \le c_{N,s} \int_{\R^N \setminus
  B_\rho(0)}|y|^{-N-2s}\,dy = \frac{N \omega_N c_{N,s}}{2s}
\rho^{-2s}=:J_{N,s}
$$
Hence  
\begin{align*}
&\int_{\cH_2}(v(t,x)-v(t,y))(\varphi(x)-\varphi(y))J(x,y)\ dxdy\\
&=2\int_{D}\varphi(x)\!\!\int_{H\setminus
  B_{\rho}(x)}\!\!(v(t,x)-v(t,y))J(x,y)\ dydx\\
& =2\int_{D}\varphi(x)\Biggl\{v(t,x)\!\!\int_{H\setminus
  B_{\rho}(x)}\!\!J(x,y) \,dydx -\!\! \int_{H\setminus
  [B_{\rho}(x) \cup \cK]  }\!\!v(t,y)J(x,y) \,dydx\\
&\qquad\qquad\qquad\qquad\qquad\qquad\qquad\qquad\qquad\qquad\qquad\qquad-\int_{\cK} v(t,y)J(x,y)\ dydx \Biggr\}\\ 
& \le 2 J_{N,s} \|\phi\|_{L^1(U_z)}\Bigl(\|v(t)\|_{L^\infty(U_z)}+\|v^-(t)\|_{L^\infty(H)}\Bigr)  - d m(t) 
\end{align*}
where in the last step we have set
$$
m(t):= \inf \limits_{y \in \cK}v(y,t)\qquad \text{and}\qquad
d:=\int_{D}\varphi(x)\int_{\cK}J(x,y)\ dydx>0.
$$
We now consider the function $t \mapsto h(t)= \int_{U_z}v(t,x)\phi(x)\,dx$ for $t>0$. Combining the
estimates above and using (\ref{eq:3}), we get 
\begin{align}
h'(t) &= \int_{\Omega_{\lambda_0}} \partial_t v(t,x)\phi(x)\,dx
\geq \int_{D}c_{\lambda_0} (t,x) v(t,x)\varphi(x)\ dx
-\mathcal{E}(v(t),\varphi)\notag\\
&\geq -c_\infty
\|\phi\|_{L^1(U_z)} \|v(t)\|_{L^\infty(U_z)}
-\mathcal{E}(v(t),\varphi) \label{eq:65}\\
&\ge -C_1 \|v(t)\|_{L^\infty(U_z)}-C_2
\|v^-(t)\|_{L^\infty(H)}  + m(t) d\nonumber
\end{align}
with $C_1:= \|\phi\|_{L^1(U_z)} \bigl[2\kappa_{s} +  c_{N,s} M  \, |D_\rho| +  J_{N,s}\bigr] + 
\frac{2|D_\rho|^2}{(2\rho)^{N+2s}} \, \|\phi\|_{L^\infty(U_z)}$ and
$C_2:=  J_{N,s} \|\phi\|_{L^1(U_z)}.$ We now consider a sequence $(t_k)_k \subset (0,\infty)$ such that $t_k
\to \infty$ and $u(t_k) \to z$ in $L^\infty(\Omega)$ as $k \to
\infty$, which yields in particular that $h(t_k) \to 0$ as $k \to \infty$. 
Using (\ref{eq:63}) and the equicontinuity property $(U2)$, we find $\delta >0$ and $k_0 \in \N$ such that
$$
m_*\,:= \inf \{m(t)\::\: t \in [t_k-\delta,t_k +\delta],\:k \ge k_0\} >0.
$$
Moreover, making $\delta>0$ smaller and $k_0 \in \N$ larger if
necessary, we may assume that 
\begin{equation}
  \label{eq:43}
\|v(t)\|_{L^\infty(U_z)} \le \|v(t)-v(t_k)\|_{L^\infty(U_z)}
+\|v(t_k)\|_{L^\infty(U_z)} \le \frac{m_*\, d}{4C_1}
\end{equation}
for $t \in [t_k-\delta,t_k +\delta]$ and $\:k \ge k_0.$ Furthermore, using that $\|v^-(t)\|_{L^\infty(H)} \to 0$ as $t \to
\infty$ as a consequence of $(S_{\lambda_0})$, we may again make $k_0
\in \N$ larger such that 
\begin{equation}
  \label{eq:64}
\|v^-(t)\|_{L^\infty(H)} \le  \frac{m_*\, d}{4C_2} \qquad \text{for $t \in [t_k-\delta,t_k +\delta],\:k \ge k_0.$}
\end{equation}
Combining (\ref{eq:65}),~(\ref{eq:43}) and (\ref{eq:64}), we thus obtain 
$$
h'(t) \ge \frac{m_*\, d}{2} \qquad \text{for $t \in [t_k-\delta,t_k
  +\delta],\:k \ge k_0.$}
$$
This implies that 
$$
\limsup_{k \to \infty} h(t_k-\delta) \le \lim_{k \to \infty}
\Bigl(h(t_k) -\frac{\delta m_*\, d}{2}\Bigr)= - \frac{\delta
  m_*\, d}{2}, 
$$
contradicting the fact that $\|v^-(t)\|_{L^\infty(U_z)} \to 0$ as $t
\to\infty$ and thus $\liminf_{t \to \infty} \limits h(t) \ge 0$. The
proof of (ii) is finished.\\
(iii) Suppose that $\lambda_0>0$, and let $z\in \omega(u)$ such that $V_{\lambda_{0}}z \equiv 0$ on $\R^N$. In view of
(ii), we need to show that $z \equiv 0$ on $\R^N$. For this we
consider the reflected functions 
\begin{align}
\tilde u&: (0,\infty) \times \R^N \to \R, && \tilde u(t,x)=
u(t,Q_0(x)) \label{eq:47}\\
\tilde z&: \R^N \to \R,&& \tilde z(x)=z(Q_0(x)).\nonumber
\end{align}
Since $\Omega$ and the nonlinearity $f$ are symmetric in the $x_1$-variable,
$\tilde u$ is also a solution of $(P)$ satisfying the same hypotheses as
$u$. Moreover, $\tilde z \in \omega(\tilde u)$. Putting $\lambda_*:=
l-2 \lambda_0 \in (-l,l)$, it follows from
$V_{\lambda_{0}}z \equiv 0$ on $\R^N$ that $\tilde z \equiv 0$ on $\Omega_{\lambda_*}$ and therefore 
\begin{equation}
    \label{eq:69}
V_\lambda \tilde z \equiv 0\quad \text{in $\Omega_{\lambda}$}\qquad \text{for every $\lambda \in (\frac{\lambda_*+l}{2},l)$.}   
  \end{equation}
For $\lambda \in (\frac{\lambda_*+l}{2},l)$ sufficiently close to $l$, it also
follows from Lemma~\ref{sec:schritt1} that $(S_\lambda)$ holds for
$\tilde u$ in place of $u$, so that (\ref{eq:69}) and (ii) imply that 
\begin{equation}
\label{eq:70}
\text{$V_\lambda \tilde z \equiv 0$ on $\R^N$ for $\lambda<l$
  sufficiently close to $l$.}   
  \end{equation}
>From this we easily conclude that $\tilde z \equiv 0$ and therefore $z
\equiv 0$ on $\R^N$, as claimed.
\qedhere
\end{proof}

\begin{lemma}\label{sec:schritt3}
 Suppose $\lambda_{0}\in (0,l)$ is such that $(S_{\lambda})$ holds for
 all $\lambda\in (\lambda_{0},l)$. Suppose furthermore that 
\underline{one} of the following conditions hold: 
\begin{itemize}
\item[(i)] $z \not \equiv 0$ on
 $\Omega$ for all $z \in \omega(u)$. 
\item[(ii)] $\Omega$ fulfills $(D2)$ and $V_{\lambda_{0}}z>0$ on
 $\Omega_{\lambda_0}$ for some $z\in \omega(u)$. 
\end{itemize}
Then there exists $\epsilon>0$ such that $(S)_{\lambda}$ holds for each $\lambda \in (\lambda_{0}-\epsilon,\lambda_{0}]$.
\end{lemma}

For the proof of this lemma, the following observation is useful. 

\begin{lemma}
\label{sec:proof-main-symmetry-1}
Let $M \subset C_0(\Omega)$ be a bounded and equicontinuous
subset, and let 
$$
I_\lambda(M):= \inf_{u \in M,\,x \in \Omega_\lambda}V_\lambda u(x) \qquad
\text{for $\lambda \in [0,l)$.}
$$
Then the map $\lambda \mapsto I_\lambda(M)$ is left continuous, i.e., for $\lambda_0
\in (0,l)$ we have $I_\lambda(M) \to I_{\lambda_0}(M)$ as $\lambda \to
  \lambda_0,\lambda< \lambda_0$.  
\end{lemma}

\begin{proof}
We first note that 
\begin{equation}
  \label{eq:72}
I_\lambda(M) \le 0 \qquad \text{for all $\lambda \in (0,l)$},  
\end{equation}
since $\overline {\Omega_\lambda} \cap T_\lambda \not=\varnothing$
by assumption $(D1)$. Since $\Omega_{\lambda_0} \subset \Omega_\lambda$ for
$\lambda<\lambda_0$ and $V_\lambda z \to V_{\lambda_0} z$ uniformly on
$\Omega_{\lambda_0}$ for every $z \in M$, we have 
$\limsup  \limits_{\lambda \to \lambda_0^-}I_\lambda(M) \le I_{\lambda_0}(M)$. Now
suppose by contradiction that there exists sequences of numbers
$\lambda_n \in (0,\lambda_0)$, of functions $u_n \in M$ and of points
$x^n \in \Omega_{\lambda_n}$ such that 
$$
\lambda_n \to \lambda \quad \text{and}\quad V_{\lambda_n}u_n(x^n) \to
c < I_{\lambda_0}(M) \qquad \text{for $n \to \infty$.} 
$$
By compactness and equicontinuity, 
we may assume that there exists $\bar x \in
\overline{\Omega}$ with ${\bar x}_1 \ge \lambda_0$ and 
$\bar u \in \overline M \subset C_0(\Omega)$ such that 
$$
x^n \to \bar x\quad  \text{and}\quad \|u_n - \bar u\|_{L^\infty(\Omega)}
\to 0 \qquad \text{as $n \to
  \infty$,}
$$
where $\overline M$ denotes the closure of $M$ in $C_0(\Omega)$ with respect to
$\|\cdot\|_{L^\infty}$. Consequently, 
$$
Q_{\lambda_n}(x^n)= 
(2\lambda_n-x^n_1,x^n_2,\dots,x^n_N) \to (2\lambda_0-\bar x_1,\bar
x_2,\dots,\bar x_N)=Q_{\lambda_0}(\bar x)
$$ 
and therefore 
$$
u_n(x^n) \to {\bar u}(\bar x) \quad \text{and}\quad u_n(Q_{\lambda_n}(x^n))
\to {\bar u}(Q_{\lambda_0}(\bar x)) \qquad \text{as $n \to \infty$.}
$$
Hence 
\begin{equation}
  \label{eq:68}
V_{\lambda_{0}}{\bar u}(\bar x)=\lim_{n \to \infty}
V_{\lambda_n}u_n(x_n)=c<I_{\lambda_0}(M) 
\end{equation}
We now distinguish two cases. If $\bar x \in \Omega$, then $\bar x \in
\overline{\Omega_{\lambda_0}}$, and we
conclude that
$$
V_{\lambda_{0}}{\bar u}(\bar x) \ge \inf_{u \in \overline M,\,x \in
\overline{\Omega_{\lambda_0}}} V_{\lambda_0} u(x) = 
\inf_{u \in M,\,x \in
\Omega_{\lambda_0}} V_{\lambda_0} u(x) = I_{\lambda_0}(M).
$$
If, on the other hand, $\bar x \in \partial
\Omega \setminus \overline{\Omega_{\lambda_0}}$, then ${\bar x}_1=
\lambda_0$ and therefore
$$
V_{\lambda_{0}}{\bar u}(\bar x)=0 \ge I_{\lambda_0}(M)
$$
by (\ref{eq:72}). Since in both cases we arrived at a statement
contradicting (\ref{eq:68}), the
proof is finished. 
\end{proof}

\begin{proof}[Proof of Lemma~\ref{sec:schritt3}]
{\em Case one:} We first assume in addition that $z \not \equiv 0$ on
$\Omega$ for all $z \in \omega(u)$. By Lemma~\ref{sec:schritt2} this
implies that $V_{\lambda_0}z>0$ in $\Omega_{\lambda_0}$ for all $z \in
\omega(u)$. Let $\delta>0$ be such that the conclusion of
Lemma~\ref{sec:schritt1} holds, and let $K \subset \Omega_{\lambda_0}$
be a compact subset and $\eps_1 \in (0,\lambda_0)$ be chosen such that 
\begin{equation}
  \label{eq:62}
|\Omega_{\lambda} \setminus
K|<\delta \qquad \text{for $\lambda \in (\lambda_0-\eps_1,\lambda_0]$.}
\end{equation}
Since $V_{\lambda_0}z>0$ in $\Omega_{\lambda_0}$ for all $z \in
\omega(u)$ and $\omega(u)$ is a compact subset of $C(\overline
\Omega)$, we may choose $\eps
\in (0,\eps_1)$ such that 
\begin{equation}
  \label{eq:60}
\inf_{z \in \omega(u) ,\,x \in K} V_\lambda z(x) >0 \qquad \text{for all $\lambda \in 
(\lambda_0-\eps,\lambda_0]$.}
\end{equation}
Let $\lambda \in (\lambda_0-\eps,\lambda_0]$, then (\ref{eq:60})
implies that there exists $t_0=t_0(\lambda)$ such that 
$$
V_{\lambda}u(t,x) \ge 0 \qquad \text{for $x \in K,\:t \ge t_0$.}
$$
Hence $\|(V_\lambda u)^-(t)\|_{L^\infty(H_\lambda)} \to 0$ as $t \to
\infty$ by Lemma~\ref{sec:schritt1}. Thus $(S_\lambda)$ holds for
$\lambda \in (\lambda_0-\eps,\lambda_0]$, as claimed.\\
{\em Case two:} We assume that $(D2)$ holds, and that $V_{\lambda_{0}}z>0$ on
 $\Omega_{\lambda_0}$ for some $z\in \omega(u)$. By $(D2)$, the set $\Omega_{\lambda_0}$ has only finitely many
connected components, and hence $\rho:=\inrad(\Omega_{\lambda_{0}})/4>0$. Let
$\gamma=\gamma(N,s,\rho,c_\infty)>0$, $q=q(N,s,\rho,c_\infty)>0$ be as in
Proposition~\ref{sec:lower-bounds-via},  and let $\delta>0$ be such that the conclusions
of Proposition \ref{sec:klein1-1} hold with $\gamma+1$  in place of $\gamma$.\\
Choose $D\subset\subset \Omega_{\lambda_{0}}$ such that $D$
intersects each connected component of $\Omega_{\lambda_0}$ and 
\begin{align}\label{sec:supi1}
 |\Omega_{\lambda_{0}}\setminus \overline{D}|<\frac{\delta}{2},\qquad
 \inrad(D)>2 \rho.
\end{align}
Fix $z\in \omega(u)$ such that $V_{\lambda_{0}}z>0$ in
$\Omega_{\lambda_{0}}$, and let $t_{n}\to \infty$ be a sequence with
$h(t_{n})\to z$. Using the equicontinuity property $(U2)$ we can find $r_{1}>0$, $\tau \in(0,\frac{1}{8})$ and $n_{0}$ such that
\begin{equation}\label{sec:supi3}
 V_{\lambda_{0}}u(t,x)>2r_{1}, \text{ for all } x\in
 \overline{D},\,t\in[t_{n}-8\tau ,t_{n}] ,\,n>n_{0}.
\end{equation}
Let $r_0:=\frac{1}{4}\dist(\overline{D},\partial
\Omega_{\lambda_{0}})$, $R= \diam(D)$  and choose $\mu$ as in Theorem
\ref{sec:haupt} for these parameter values. We first fix $\eps_1>0$
such that  
\begin{equation}\label{sec:supi2}
 |\Omega_{\lambda}\setminus\Omega_{\lambda_{0}}|<\frac{\delta}{2}, \text{ for } \lambda\in[\lambda_{0}-\epsilon_1,\lambda_{0}).
\end{equation}
>From the equicontinuity assumption (U2) we may deduce that 
\begin{equation}
\sup_{n \in \N} \sup_{[t_{n}-8 \tau,t_{n}]\times D}|V_{\lambda}u-V_{\lambda_{0}}u|\to 0 \text{ as } \lambda\to \lambda_{0}.
\end{equation}
This and (\ref{sec:supi3}) imply the existence of
$\eps_2 \in (0,\eps_1)$ such that
\begin{equation}\label{sec:supi4}
 V_{\lambda}u(t)>r_{1}, \text{ for all } x\in
 \overline{D},\,t\in[t_{n}-8 \tau,t_{n}] ,\,n>n_{0}, \, \lambda\in [\lambda_{0}-\epsilon_2,\lambda_{0}].
\end{equation}
By $(S_{\lambda_{0}})$, we can find $n_{1}>n_{0}$ such that for all $n>n_{1}$ we have
\[
 \|(V_{\lambda_{0}}u)^{-}(t_{n}-8\tau)\|_{L^{\infty}(\Omega_{\lambda_{0}}\setminus \overline{D})}\leq \frac{\mu r_{1}}{2}.
\]
Using the equicontinuity of the functions $x \mapsto u(t_n-8
\tau,x)$, $n \in \N$ and Lemma~\ref{sec:proof-main-symmetry-1}, we may
choose $\eps \in (0,\eps_2)$ such that
\begin{equation}
  \label{eq:41}
\|(V_{\lambda}u)^{-}(t_{n}-8\tau)\|_{L^{\infty}(\Omega_{\lambda}\setminus
   \overline{D})}\leq \mu r_{1}\qquad \text{for $\lambda\in [\lambda_{0}-\epsilon,\lambda_{0}].$} 
\end{equation}
 We now fix $n \ge n_1$ and $\lambda\in
[\lambda_{0}-\epsilon,\lambda_{0}]$, and we claim that the assumptions of
Theorem~\ref{sec:haupt} are satisfied with $t_0=t_n-8\tau$, $U=\Omega_\lambda$, $D$ as
above and $v=
V_\lambda u$. Indeed, $\dist(\overline{D},\partial U)\geq
\dist(\overline{D},\partial\Omega_{\lambda_{0}}) \ge 4r_0$ and $|\Omega_\lambda \setminus
D|<\delta$ by (\ref{sec:supi1}) and (\ref{sec:supi2}). Moreover,
$\inrad(D) >2\rho$ and $\diam D \le R$ by our choice of $D$ and the
definition of $R$. Moreover, by (\ref{sec:supi4}), $V_\lambda u$ is nonnegative on
$[t_n-8 \tau ,t_{n}] \times \overline{D}$, and by (\ref{sec:supi4})
and (\ref{eq:41}) we
have 
$$
\|(V_{\lambda}u)^{-}(t_n-8 \tau)\|_{L^{\infty}(U\setminus
   \overline{D})} \leq \mu r_{1} \le \mu
 [V_{\lambda}u]_{L^1([t_n-7\tau,t_n-6\tau] \times D_*)}.  
$$
for each connected component $D_*$ of $D$. An application of Theorem \ref{sec:haupt}(ii) with these parameters
therefore yields that $(S_{\lambda})$ holds for all $\lambda\in
[\lambda_{0}-\epsilon,\lambda_{0}]$. The proof is finished.
\end{proof}

The following Proposition evidently completes the {\em Proof of Theorem \ref{sec:goal}}.

\begin{proposition}
\label{sec:proof-main-symmetry-2}
Suppose that $(D2)$ holds or that $z \not \equiv 0$ on
 $\Omega$ for all $z \in \omega(u)$. Then we have:
\begin{itemize}
\item[(i)] $V_{0}z \equiv 0$ on $\R^N$ for
every $z \in \omega(u)$.
\item[(ii)] For every $z \in \omega(u)$, we either have the following
  alternative. Either $z \equiv 0$ on $\Omega$,
  or $z$ is strictly decreasing in $|x_1|$ and therefore strictly
  positive in $\Omega$.
\end{itemize}
\end{proposition}

\begin{proof}
(i) We define  
$$\lambda_{0}:=\inf\{\mu >0: (S_{\lambda}) \text{ holds for all }
\lambda>\mu\},
$$
and we first claim that $\lambda_0=0$. By Lemma~\ref{sec:schritt1} we
have $\lambda_0<l$. If $z \not \equiv 0$ on
 $\Omega$ for all $z \in \omega(u)$, then Lemma~\ref{sec:schritt3}
 immediately implies that $\lambda_0=0$. If $(D2)$ holds and we assume
 -- on the contrary --  $\lambda_0>0$, then
Lemma~\ref{sec:schritt2}(iii) and Lemma~\ref{sec:schritt3}(ii) readily
imply that $z \equiv 0$ on $\R^N$ for
every $z \in \omega(u)$, which then also yields $\lambda_0=0$. Hence
we conclude in both cases that
$\lambda_0=0$, and therefore $(S_0)_0$ is true by
Lemma~\ref{sec:schritt2}(i). This implies that $V_0 z
\ge 0$ on $\Omega_0$ for every $z \in \omega(u)$. Since the analogous
statement can also be shown for the reflected
solution $\tilde u$ defined in (\ref{eq:47}), we also
have that $V_0 z \le 0$ on $\Omega_0$ for every $z \in
\omega(u)$. Hence for every $z \in \omega(u)$ we have $V_0 z \equiv 0$
on $\Omega_0$ and thus also on $\R^N$, since $z \equiv 0$ on $\R^N
\setminus \Omega$.\\
(iii) Let $z \in \omega(u)$ be given such that $z$ is not strictly
decreasing in $|x_1|$. Then there exists $\lambda>0$ such that 
$V_\lambda z$ is not strictly positive in $\Omega_\lambda$. By
Lemma~\ref{sec:schritt2}(ii), applied to $\lambda$ in place of
$\lambda_0$, we then have that $V_\lambda z \equiv 0$ on
$\R^N$. By (ii), $z$ therefore has two different parallel symmetry
hyperplanes. This implies that $z \equiv 0$, since $z$ vanishes
outside a bounded subset of $\R^N$.
\end{proof}
 
\section{Appendix}
\label{sec:appendix}

%Änderung Sven Jarohs 9.11.

As announced in the introduction, we derive -- based on recent results in \cite{FK12}
and \cite{RS12} -- a sufficient
criterion for condition $(U2)$. For a similar result in the context of local parabolic
boundary value problems, see \cite[Prop. 2.7]{P07}.

\begin{proposition}
\label{sec:appendix-6}
 Let $\Omega\subset\R^{N}$ be a bounded domain, and suppose that the
 nonlinearity $f$ satisfies $(F1)$. Suppose furthermore that $0\in \mathcal{B}$, and that $f(\cdot,\cdot,0)$
 is bounded on $(0,\infty)\times \Omega$. Then for any solution $u$
 of $(P)$ satisfying $(U1)$ we have: 
 \begin{itemize}
 \item[(i)]  For any domain
 $G\subset \subset\Omega$ there exist $\alpha>0$ such that
 \begin{equation}
   \label{eq:58}
  \sup_{\substack{\tau \ge 1\\ t,\tilde{t}\in [\tau,\tau+1], t\neq
      \tilde{t}\\  x,\tilde{x}\in\overline{G},x\neq \tilde{x}}}\frac{|u(t,x)-u(\tilde{t},\tilde{x})|}{\left(|x-\tilde{x}|+|t-\tilde{t}|^{1/2s}\right)^\alpha}<\infty.
 \end{equation}
 \item[(ii)] If, in addition, $\Omega$ fulfills the exterior sphere
   condition and, for some $t_0 > 0$, $C_1 >0$, 
   \begin{equation}
     \label{eq:51}
|u(t_0,x)| \le C_1 \dist(x,\partial \Omega)^s \qquad \text{for all $x
  \in \Omega$,}     
   \end{equation}
then 
   \begin{equation}
\label{eq:52}
\sup_{t \ge t_0,\, x \in \Omega}\; \frac{|u(t,x)|}{\dist(x,\partial
  \Omega)^s} <\infty
   \end{equation}
In particular, $(U2)$ holds in this case.
 \end{itemize}
\end{proposition}

In the special case $f \equiv 0$, the interior 
regularity estimate (\ref{eq:58}) is an immediate consequence of
\cite[Theorem 1.2]{FK12}, but we could not find any reference where
the case $f \not \equiv 0$ is considered. Before giving the proof of this proposition, we discuss an example. 

\begin{remark}
Let $\Omega \subset \R^N$ be a bounded
domain satisfying the exterior sphere condition. We consider an
Allen-Cahn-type nonlinearity 
\begin{equation}
  \label{eq:50}
f:[0,\infty) \times \Omega \times \R \to \R, \qquad f(t,x,u)= a(t)u
-b(t) u^3= u[a(t)-b(t)u^2]
\end{equation}
Here $a,b: [0,\infty) \to \R$ are continuous functions with $a(t) \le
b(t)$ for $t \ge 0$. Then $f$ satisfies $(F1)$ with $\cB=\R$, and it trivially 
satisfies $(F2)$ if
$\Omega$ satisfies $(D1)$. Moreover, the constant $1$ is a supersolution of   
problem $(P)$, whereas $0$ is a solution. Hence, if $\phi \in
C_0(\Omega) \cap \mathcal{H}^{s}_{0}(\Omega)$ is such that $0 \le \phi(x) \le 1$ for all $x \in
\Omega$, standard methods in semigroup theory and the weak maximum
principle (see Remark~\ref{sec:small-volume-maximum}) 
give rise to the existence of a unique global solution of the initial value problem 
\begin{equation}
  \label{eq:44}
\left\{\begin{aligned}
& u\in C([0,\infty),\mathcal{H}^{s}_{0}(\Omega) \cap C_{0}(\Omega))\cap
C^{1}((0,\infty),L^{2}(\Omega)),\quad  \hl u \in C((0,\infty),L^{2}(\Omega))\\
      &\partial_{t}u(t)+\hl u(t)= f(t,x,u(t)) \qquad\qquad \text{ for $t\in (0,\infty)$,}\\
&u(0)=\varphi.
     \end{aligned}  
\right.
\end{equation}
satisfying $0 \le u(t,x) \le 1$ for all $t \in (0,\infty)$, $x \in
\Omega$, so that condition $(U1)$ is satisfied for $u$. Furthermore,
if $\phi(x) \le C_1 \dist(x,\partial \Omega)^s$ for $x
  \in \Omega$ with some constant $C_1>0$, then $(U2)$ is also satisfied by Proposition~\ref{sec:appendix-6}(ii).
We remark that the solution $u$ can be found as a the unique mild
solution of (\ref{eq:44}), i.e.,
the unique solution of the nonlinear integral equation
\begin{equation}
  \label{eq:30}
u \in C([0,\infty),C_{0}(\Omega)),\qquad 
u(t)= S_{A}(t)\varphi + \int_{0}^{t}S_{A}(t-\tau)F(\tau,u(\tau))\
d\tau \quad \text{for $t\in[0,\infty)$}.
\end{equation}
Here $S_A$ denotes the semigroup generated by the $m$-dissipative operator 
$$
A: \dom(A) \subset C_{0}(\Omega) \to C_{0}(\Omega),\qquad Au:=-\hl u 
$$
where $\dom(A)$ is the space of all functions $u\in\mathcal{H}_{0}^{s}(\Omega)\cap
C_{0}(\Omega)$ such that $\hl u$, defined in distributional sense, is
contained in $C_{0}(\Omega)$. Moreover, $F:[0,\infty) \times C_0(\Omega) \to
C_0(\Omega)$ is the substitution operator given by $[F(t,w)](x)=
f(t,x,w(x))$ for $t \in [0,\infty)$, $x \in \Omega$. The
$m$-dissipativity of the operator $A$ in $C_0(\Omega)$ is essentially
a consequence of the following recent regularity result given in \cite[Proposition
1.1]{RS12}: If $\Omega \subset \R^N$ is a bounded domain satisfying
the exterior sphere condition and $w \in L^\infty(\Omega)$, then the unique weak solution $u
\in H^s_0(\Omega)$ of the equation $-\Delta u=w$ belongs to
$C_0(\Omega)$. Another important fact needed for the local existence
and uniqueness of solutions of $(P)$ is the local uniform (in time) Lipschitz continuity of $F:[0,\infty) \times C_0(\Omega) \to
C_0(\Omega)$, which follows since $f$ satisfies $(F2)$. In
order to show solutions of (\ref{eq:30}) are also
solutions of (\ref{eq:44}), one may essentially argue as in
\cite{CH98} for the semilinear heat equation, noting the following
additional useful property of the substitution operator $F$: If $M \subset C_0(\Omega) \cap
\mathcal{H}^{s}_{0}(\Omega)$ is bounded with respect to
$\|\cdot\|_\infty$, then $F(M) \subset \mathcal{H}^{s}_{0}(\Omega)$,
and there exists $L=L(M)>0$ such that 
\begin{equation}
  \label{eq:48}
\cE(F(t,u),F(t,u)) \le L \cE(u,u) \qquad \text{for all $u \in M$,
  $t>0$.}  
\end{equation} 
This property can be checked immediately by using $(F2)$ and the definition of
the quadratic form $\cE$. 

Note that (\ref{eq:50}) is
just a particular example of a nonlinearity which admits an ordered
pair of a bounded subsolution $\phi_*$ and a bounded supersolution
$\phi^*$ and which
satisfies $(F1)$ with $\cB=\R$. In such a setting, an initial
condition $\phi \in C_0(\Omega) \cap \mathcal{H}^{s}_{0}(\Omega)$
always gives rise to a global bounded solution of $(P)$.
\end{remark}

The remainder of this appendix is devoted to the proof of
Proposition~\ref{sec:appendix-6}. The assertion (\ref{eq:58}) on
interior regularity will be deduced from the Harnack inequality of
Felsinger and Kassmann \cite{FK12}. More precisely, we will use the
following rescaled variant of a special case of \cite[Corollary
5.2]{FK12}. 

\begin{proposition}
\label{sec:appendix-1}
Let 
$$
D_{\ominus} := (-2^{2s+1},-2^{2s+1}+1) \times
B_{1}(0)\qquad \text{and}\qquad D_{\oplus} := (-1,0) \times
B_{1}(0).
$$
There exists $\eps_0,\delta>0$ such that for every
nonnegative supersolution 
$$
w:(-2^{2s+1},0) \times \R^N \to \R
$$
of the equation 
$$
\partial_t w+(-\Delta)^s = -\eps_0 \qquad \text{in $(-2^{2s+1},0) \times
  B_4(0)$}
$$
in the sense of Definition~\ref{sec:antisymm-supers-corr-2} with the
property that 
\begin{equation}
  \label{eq:73}
|D_{\ominus} \cap \{w \ge 1\}| \ge \frac{1}{2} |D_{\ominus}|   
\end{equation}
we have $w \ge \delta$ a.e. on $D_{\oplus}$.
\end{proposition}

\begin{corollary}
\label{hoelder-reg-1}
 Let $r_{0} \in (0,1]$, $c_{u}>0$ and
 $f_{\infty}> 0$. Then there exist constants $\alpha \in
 (0,1)$ and $C_{2}>0$ depending on $N,s,f_{\infty},c_{u},r_{0}$ with the
 following property:\\
If $T:=(t_0-r_0^{2s},t_0)$ for some $t_0 \in \R$, $x_0 \in \R^N$, $f \in
L^\infty(T\times B_{r_0}(x_0))$ with $\|f\|_{L^\infty(T\times
  B_{r_0}(x_0))}\le f_\infty$ are given and 
$$
u \in
C(T,H^{s}(\R^N)\cap L^\infty(\R^N) \cap C(\overline{B_{r_0}(x_0)}) \cap C^1(T,L^2(B_{r_0}(x_0)))   
$$
with $\|u\|_{L^\infty (T \times \R^N)}\le c_u$ is a solution of 
$$
\partial_t u + \hl
u= f(t,x)\qquad \text{in $T\times B_{r_0}(x_{0})$}
$$ 
in the sense that 
$$
\cE(u(t),\varphi)=
\int_{B_{r_0}(x_0)}\!\!\left[f(t,x)-\partial_{t}u(t,x)\right]\varphi(x)\
dx
$$
for every $\varphi \in \cH^{s}_{0}(B_{r_0}(x_0))$ and a.e. $t \in T$, then we have 
\begin{equation}
  \label{eq:49}
\underset{Q(r)}{\osc}u \le C_2 r^\alpha \qquad \text{for $r \in
  (0,r_0]$},\; \text{where $Q(r):=(t_0- r^{2s},t_0) \times B_{r}(x_0)$.}  
\end{equation}
\end{corollary}

\begin{proof}
Without loss, we may assume that $t_0=0$ and $x_0 =0$. Moreover, we
may assume by normalization that $c_u = \frac{1}{4}$. In this case we
will prove (\ref{eq:49}) with $C_2=1$ for some suitable $\alpha \in
(0,1)$. Suppose by contradiction that the statement is false. Then there exist,
for every $k \in \N$, functions $f_k \in
L^\infty(T\times B_{r_0}(0))$ with $\|f_k\|_{L^\infty(T\times
  B_{r_0}(0))}\le f_\infty$ and $u_k \in
C(T,H^{s}(\R^N)\cap L^\infty(\R^N) \cap C(\overline{B_{r_0}(0)}) \cap C^1(T,L^2(B_{r_0}(0)))   
$
with 
$$
\|u_k\|_{L^\infty(T \times \R^N)}\le \frac{1}{4}
$$
solving 
$$
\partial_{t}u_k + \hl
u_k= f_k(t,x)\qquad \text{in $T \times B_{r_0}(0)$}
$$
as well as $\alpha_k \in (0,1)$ and $r_k \in (0,r_0]$ 
such that $\alpha_k \to 0$ as $k \to \infty$ and
$$
\underset{Q(r_k)}{\osc}u_k \ge r_k^{\alpha_k} \qquad \text{for every $k
  \in \N$.} 
$$
Passing to a subsequence, we also have   
$$
\underset{T \times \R^N}{\osc}\,u_k \le 2 \|u_k\|_{L^\infty(T
  \times \R^N)} \le \frac{1}{2} \le r_0^{\alpha_k} \qquad \text{for every $k \in \N$.}
$$
By making $r_k \in (0,r_0]$ larger if necessary, we may therefore
assume that  
$$
\underset{Q(r_k)}{\osc}u_k = r_k^{\alpha_k} \qquad \text{for every $k
  \in \N$} 
$$
and 
$$
\underset{Q(r)}{\osc}\,u_k \le r^{\alpha_k} \qquad \text{for $r \in [r_k,r_0]$
  and $k \in \N$.} 
$$
Since also $\underset{Q(r_k)}{\osc}\,u_k \le
\frac{1}{2}$ for every $k \in \N$, we conclude that $r_k \to 0$ as $k \to \infty$. We now
define $T_k:=(-(\frac{r_0}{r_k})^{2s},0)$ and 
$$
v_k: T_k \times \R^N \to \R, \qquad v_k(t,x)= 2r_k^{-\alpha_k} u_k(r_k^{2s}t,r_k
x)
$$
for $k \in \N$. Then we have
$$
\partial_{t}v_k + \hl
v_k= \tilde f_k(t,x)\qquad \text{in $D_k:=T_k \times B_{\frac{r_0}{r_k}}(0)$}
$$
with 
$$
\tilde f_k(t,x)= 2r_k^{2s-\alpha_k}f_k(r_k^{2s}t,r_k x).
$$
Without loss, we may assume that $\frac{r_0}{r_k} \ge \max
\{2^{1+\frac{1}{2s}},5\}$ for every $k \in \N$, so that $(-2^{2s+1},0) \times B_5(0) \subset
D_k$ for every $k \in \N$. Moreover, we have $\underset{Q(1)}{\osc}\,v_k = 2$, 
\begin{equation}
  \label{eq:74}
\underset{Q(r)}{\osc}\,v_k \le 2r^{\alpha_k} \qquad \text{for $r \in [1,\frac{r_0}{r_k}]$, $k \in \N$} 
\end{equation}
and 
\begin{equation}
\label{eq:67}
\underset{T_k \times \R^N}{\osc}\,v_k \le 2 \Bigl(\frac{r_0}{r_k}\Bigr)^{\alpha_k} \qquad \text{for $k \in \N$.} 
\end{equation}
By adding a constant to $v_k$ if necessary, we may assume that
\begin{equation}
  \label{eq:75}
\underset{{Q(1)}}{\sup}\,v_k= 1 \qquad \text{and}\qquad \underset{{Q(1)}}{\inf}\,v_k=-1.
\end{equation}
After passing to a subsequence, we may also assume that, 
replacing $v_k$ by $-v_k$ and $\tilde f_k$ by $-\tilde f_k$ if necessary,
$$
|D_{\ominus} \cap \{v_k \ge 0\}| \ge \frac{1}{2} |D_{\ominus}|.   
$$
Here and in the following, $D_\ominus$ and $D_\oplus$ are defined as in
Proposition~\ref{sec:appendix-1}. Note that by (\ref{eq:74}),~(\ref{eq:67}) and (\ref{eq:75}) we have
$$
v_k(t,x) \ge \min \{-1, 1-2|x|^{\alpha_k}\} \qquad \text{for $x \in
  \R^N$, $\:t \in (-2^{2s+1},0)$.} 
$$
We now consider 
$$
w_k: T_k \times \R^N \to \R,\qquad w_k(t,x):=
v_k(t,x)+ 2 \cdot 5^{\alpha_k}-1.
$$
Then 
$$
w_k(t,x) \ge \min \bigl\{0, 2(5^{\alpha_k}-
|x|^{\alpha_k})\bigr\} \qquad \text{for $x \in \R^N$, $\:t \in (-2^{2s+1},0)$.} 
$$
In particular, we have $w_k \ge 0$ in $(-2^{2s+1},0) \times B_5(0)$, and for
$x \in B_4(0)$ we have 
$$
|(-\Delta)^s w_k^-(t,x)|\le 2 \int_{\R^N \setminus
  B_5(0)}\frac{|y|^{\alpha_k}-5^{\alpha_k}}{|x-y|^{N+2s}}\,dy
\le \int_{\R^N \setminus
  B_5(0)}\frac{|y|^{\alpha_k}-5^{\alpha_k}}{(|y|-4)^{N+2s}}\,dy, 
$$
where the latter integral tends to zero as $k \to \infty$ by
Lebesgue's theorem. Hence  
\begin{equation}
  \label{eq:71}
\lim_{k \to \infty}\|(-\Delta)^s w_k^-\|_{L^\infty((-2^{2s+1},0)\times
  B_4(0))} = 0.
\end{equation}
We now note that the function $w_k^{+}$ is a nonnegative solution of 
\[
 \partial_{t}w_k^{+}+\hl w_{k}^{+}=g_k \qquad\text{in $(-2^{2s+1},0)
   \times B_4(0)$ for every $k\in \N$}
\]
with $g_k:= \tilde{f}_{k}+\hl w_{k}^{-}$, whereas $\|g_k\|_{L^\infty((-2^{2s+1},0)
   \times B_4(0))} \to 0$ as $k \to \infty$ as a consequence of
 (\ref{eq:71}) and the fact that 
$$
\|\tilde f_k\|_{L^\infty((-2^{2s+1},0)
   \times B_4(0))}\le 2r_k^{2s-\alpha_k} f_\infty
$$
Consequently, there exists $k_{0}\in\N$ such that $\|g_k\|_{L^\infty((-2^{2s+1},0)
   \times B_4(0))} \le \eps_0$, where $\epsilon_{0}$ is given by Lemma
 \ref{sec:appendix-1}. On the other hand, since $D_{\oplus}=Q(1)$, we
 infer from (\ref{eq:75}) that
$$
\inf_{D_\oplus} w_k^+ = \inf_{D_\oplus} w_k = 2 \cdot 5^{\alpha_k}-2 \to 0
\qquad \text{as $k \to \infty$.}
$$
This contradicts Proposition~\ref{sec:appendix-1}, applied to $w=w_k^+$. The proof is thus finished.
\end{proof}

\begin{proof}[Proof of Proposition~\ref{sec:appendix-6} (completed)]

(i) We note that $u$ satisfies
$$
\partial_{t}u(t,x)+\hl u(t,x)=\tilde f(t,x)
$$
with $\tilde f(t,x)= f(t,x,u(t,x))$, and by assumption 
$u$ and $f$ are bounded on \mbox{$(0,\infty) \times \R^N$}. Hence, for given $G
\subset \subset \Omega$, we may choose $r_0 >0$ such that $r_0 < \min
\{\dist(G,\partial
  \Omega),1\}$, and we
may apply Corollary~\ref{hoelder-reg-1} to every point $x_0 \in G$, $t_0 \ge 1$. From this
(\ref{eq:58}) easily follows.\\
(ii) We use barrier functions as constructed in the elliptic setting
in \cite{RS12}. Put $B_r:= B_r(0)$ for $r>0$, and recall the
definition of the space $\cV^s(U')$ in (\ref{eq:66}). By \cite[Lemma 2.6]{RS12} there exists
a function $\phi \in \cV^s(\R^N)$ satisfying 
\begin{equation}
  \label{eq:53}
\left \{
  \begin{aligned}
 &\hl \phi \ge 1 \quad \text{in $B_4 \setminus B_1$,}&&\qquad 0 \le
 \phi(x) \le c_0 (|x|-1)^s \quad \text{for $x \in B_4 \setminus B_1$;}\\
 &\phi \equiv 0\quad \text{in $B_1$,}&&\qquad \phi \ge 1\quad \text{in $\R^N \setminus B_4$}
 \end{aligned}
\right.
\end{equation}
as well as 
\begin{equation}
  \label{eq:54}
\phi(x) \ge d_0\: \dist(x,\partial B_1)^s \qquad \text{for $x \in B_4 \setminus
  B_1$}   
\end{equation}
with some constants $c_0, d_0>0$. In fact, it is not stated explicitly
in \cite{RS12} that $\phi \in \cV^s(\R^N)$ and
that (\ref{eq:54}) holds, but this follows from the
construction in \cite[Appendix]{RS12}. Now since $\Omega$ satisfies
the exterior sphere condition, there exists $\rho>0$ such that every
point in $\partial \Omega$ can be touched from outside by a ball of
radius $\rho$. Fixing such a ball $B_\rho(y)$ for some $y \in \R^N
\setminus \Omega$, we may define the function 
$$
\psi \in H^s_{loc}(\R^N),\qquad \psi(x)= \lambda \phi (\frac{x-y}{\rho}).
$$
Here, using (\ref{eq:51}), (\ref{eq:53}), (\ref{eq:54}) and the assumption that $u$
satisfies $(U1)$, we may choose $\lambda>0$
sufficiently large so that 
\begin{equation}
  \label{eq:55}
\left \{
  \begin{aligned}
\hl \psi &\ge\;  \sup_{t \ge t_0, x \in \Omega}f(t,x,u(t,x))  &&\qquad \text{in $B_{4\rho}(y) \setminus B_{\rho}(y)$,}\\
 \psi &\ge \;\sup_{t \ge t_0, x \in \Omega}u(t,x)  && \qquad \text{in $\R^N \setminus B_{4\rho}(y)$}\\
\psi(x) &\ge\:  u(t_0,x) && \qquad \text{for $x \in \Omega \cap B_{4\rho}(y)$.}
 \end{aligned}
\right.
\end{equation}
Let $w(t,x)= \psi(x)-u(t,x)$. By the
properties (\ref{eq:55}), $w$ is an entire supersolution of $\partial_t w + \hl w = 0$ in
$[t_0,\infty) \times [\Omega \cap B_{4\rho}(y)]$ in the sense of
Definition~\ref{sec:antisymm-supers-corr-2}, and $w(t_0)$ is
nonnegative on $\R^N$. Hence, by the weak maximum
principle as stated in Remark~\ref{sec:small-volume-maximum}, $w(t,x) \ge 0$ for $x \in \Omega$, $t \ge t_0$ and
therefore 
$$
u(t,x) \le \psi(x) \le \frac{\lambda c_0}{\rho^s} (|x-y|-\rho)^s
\qquad \text{for $x \in \Omega \cap B_{4\rho}(y), \; t\ge t_0.$}
$$
Since the parameter $\lambda$ in the definition of $\phi$ can be chosen
uniformly with respect to the $\rho$-balls touching $\Omega$ from
outside, we find -- using also the boundedness of $u$ on
$[t_0,\infty)\times \Omega$ -- a constant $C'>0$ such that 
\begin{equation}
  \label{eq:56}
u(t,x) \le C' \dist(x,\partial \Omega)^s \qquad \text{for $x \in
  \Omega$, $t \ge t_0$.}  
\end{equation}
Repeating the same argument with $-u$ in place of $u$, we find a
constant $C''>0$ such that    
\begin{equation}
\label{eq:57}
u(t,x) \ge -C'' \dist(x,\partial \Omega)^s \qquad \text{for $x \in
  \Omega$, $t \ge t_0$.}  
\end{equation}
Combining (\ref{eq:56}) and (\ref{eq:57}), we obtain (\ref{eq:52}), as
claimed. Now $(U2)$ follows easily by combining (\ref{eq:58}) and 
(\ref{eq:52}).
\end{proof}

%For acknowledgements section, please don't number the section, please begin it with \section*{Acknowledgements}
\section*{Acknowledgments} The authors would like to thank Mouhamed
Moustapha Fall and Peter Pol\'a\v cik for helpful discussions.

% You may incorporate your references as follows in your main tex file.
% Using BibTex is not recommended but can be handled.

% \medskip
% % The data information below will be filled by AIMS editorial staff
% Received xxxx 20xx; revised xxxx 20xx.
% \medskip

\end{document}